\let\noarrow = t
\input eplain

\let\noarrow = t

\input eplain


\magnification=\magstep1

\topskip1truecm
\def\pagewidth#1{
  \hsize=#1
}

\def\pageheight#1{
  \vsize=#1
}

\pageheight{23.5truecm} \pagewidth{16truecm}

\abovedisplayskip=3mm \belowdisplayskip=3mm
\abovedisplayshortskip=0mm \belowdisplayshortskip=2mm
\parindent1pc

\normalbaselineskip=13pt \baselineskip=13pt

\def\spacing{{\smallskip}}

\voffset=0pc \hoffset=0pc


\newdimen\abstractmargin
\abstractmargin=3pc


\newdimen\footnotemargin
\footnotemargin=1pc


\font\eightrm=cmr8 \relax 
\font\sixrm=cmr6 \relax 
\font\eighti=cmmi8 \relax     \skewchar\eighti='177 
\font\sixi=cmmi6 \relax       \skewchar\sixi='177   
\font\eightsy=cmsy8 \relax    \skewchar\eightsy='60 
\font\sixsy=cmsy6 \relax      \skewchar\sixsy='60   
\font\eightbf=cmbx8 \relax 
\font\sixbf=cmbx6 \relax   
\font\eightit=cmti8 \relax 
\font\eightsl=cmsl8 \relax 
\font\eighttt=cmtt8 \relax 

\catcode`\@=11
\newskip\ttglue

\def\eightpoint{\def\rm{\fam0\eightrm}%
 \textfont0=\eightrm \scriptfont0=\sixrm
 \scriptscriptfont0=\fiverm
 \textfont1=\eighti \scriptfont1=\sixi
 \scriptscriptfont0=\fivei
 \textfont2=\eightsy \scriptfont2=\sixsy
 \scriptscriptfont2=\fivesy
 \textfont3=\tenex \scriptfont3=\tenex
 \scriptscriptfont3=\tenex
 \textfont\itfam\eightit \def\it{\fam\itfam\eightit}%
 \textfont\slfam\eightsl \def\sl{\fam\slfam\eightsl}%
 \textfont\ttfam\eighttt \def\tt{\fam\ttfam\eighttt}%
 \textfont\bffam\eightbf \scriptfont\bffam\sixbf
   \scriptscriptfont\bffam\fivebf \def\bf{\fam\bffam\eightit}%
 \tt \ttglue=.5em plus.25em minus.15em
 \normalbaselineskip=9pt
 \setbox\strutbox\hbox{\vrule height7pt depth3pt width0pt}%
 \let\sc=\sixrm \let\big=\eifgtbig \normalbaselines\rm}


 \font\titlefont=cmbx12 scaled\magstep1
 \font\sectionfont=cmbx12
 \font\ssectionfont=cmsl10
 \font\claimfont=cmsl10

 \font\normalfont=cmr10

\catcode`\@=11 \font\teneusm=eusm10 
\font\seveneusm=eusm7  \font\fiveeusm=eusm5
\newfam\eusmfam \textfont\eusmfam=\teneusm
\scriptfont\eusmfam=\seveneusm \scriptscriptfont\eusmfam=\fiveeusm
\def\hexnumber@#1{\ifcase#1
0\or1\or2\or3\or4\or5\or6\or7\or8\or9\or         A\or B\or C\or D\or
E\or F\fi } \edef\eusm@{\hexnumber@\eusmfam}
\def\euscr{\ifmmode\let\next\euscr@\else
\def\next{\errmessage{Use \string\euscr\space only in math mode}}\fi\next}
\def\euscr@#1{{\euscr@@{#1}}} \def\euscr@@#1{\fam\eusmfam#1} \catcode`\@=12

\catcode`\@=11 \font\teneuex=euex10 
 \font\seveneuex=euex7  \newfam\euexfam
\textfont\euexfam=\teneuex  \scriptfont\euexfam=\seveneuex
 \def\hexnumber@#1{\ifcase#1
0\or1\or2\or3\or4\or5\or6\or7\or8\or9\or         A\or B\or C\or D\or
E\or F\fi } \edef\euex@{\hexnumber@\euexfam}
\def\euscrex{\ifmmode\let\next\euscrex@\else
\def\next{\errmessage{Use \string\euscrex\space only in math mode}}\fi\next}
\def\euscrex@#1{{\euscrex@@{#1}}} \def\euscrex@@#1{\fam\euexfam#1}
\catcode`\@=12

\catcode`\@=11 \font\teneufb=eufb10 
\font\seveneufb=eufb7  \font\fiveeufb=eufb5
\newfam\eufbfam \textfont\eufbfam=\teneufb
\scriptfont\eufbfam=\seveneufb \scriptscriptfont\eufbfam=\fiveeufb
\def\hexnumber@#1{\ifcase#1
0\or1\or2\or3\or4\or5\or6\or7\or8\or9\or         A\or B\or C\or D\or
E\or F\fi } \edef\eufb@{\hexnumber@\eufbfam}
\def\euscrfb{\ifmmode\let\next\euscrfb@\else
\def\next{\errmessage{Use \string\euscrfb\space only in math mode}}\fi\next}
\def\euscrfb@#1{{\euscrfb@@{#1}}} \def\euscrfb@@#1{\fam\eufbfam#1}
\catcode`\@=12

\catcode`\@=11 \font\teneufm=eufm10 
\font\seveneufm=eufm7  \font\fiveeufm=eufm5
\newfam\eufmfam \textfont\eufmfam=\teneufm
\scriptfont\eufmfam=\seveneufm \scriptscriptfont\eufmfam=\fiveeufm
\def\hexnumber@#1{\ifcase#1
0\or1\or2\or3\or4\or5\or6\or7\or8\or9\or         A\or B\or C\or D\or
E\or F\fi } \edef\eufm@{\hexnumber@\eufmfam}
\def\euscrfm{\ifmmode\let\next\euscrfm@\else
\def\next{\errmessage{Use \string\euscrfm\space only in math mode}}\fi\next}
\def\euscrfm@#1{{\euscrfm@@{#1}}} \def\euscrfm@@#1{\fam\eufmfam#1}
\catcode`\@=12

\catcode`\@=11 \font\teneusb=eusb10 
\font\seveneusb=eusb7  \font\fiveeusb=eusb5
\newfam\eusbfam \textfont\eusbfam=\teneusb
\scriptfont\eusbfam=\seveneusb \scriptscriptfont\eusbfam=\fiveeusb
\def\hexnumber@#1{\ifcase#1
0\or1\or2\or3\or4\or5\or6\or7\or8\or9\or         A\or B\or C\or D\or
E\or F\fi } \edef\eusb@{\hexnumber@\eusbfam}
\def\euscrsb{\ifmmode\let\next\euscrsb@\else
\def\next{\errmessage{Use \string\euscrsb\space only in math mode}}\fi\next}
\def\euscrsb@#1{{\euscrsb@@{#1}}} \def\euscrsb@@#1{\fam\eusbfam#1}
\catcode`\@=12

\catcode`\@=11 \font\tenmsa=msam10 
\font\sevenmsa=msam7  \font\fivemsa=msam5
\font\tenmsb=msbm10  \font\sevenmsb=msbm7
 \font\fivemsb=msbm5 \newfam\msafam
\newfam\msbfam \textfont\msafam=\tenmsa
\scriptfont\msafam=\sevenmsa
  \scriptscriptfont\msafam=\fivemsa
\textfont\msbfam=\tenmsb  \scriptfont\msbfam=\sevenmsb
  \scriptscriptfont\msbfam=\fivemsb
\def\hexnumber@#1{\ifcase#1 0\or1\or2\or3\or4\or5\or6\or7\or8\or9\or
        A\or B\or C\or D\or E\or F\fi }
\edef\msa@{\hexnumber@\msafam} \edef\msb@{\hexnumber@\msbfam}
\mathchardef\square="0\msa@03 \mathchardef\subsetneq="3\msb@28
\mathchardef\supsetneq="3\msb@29 \mathchardef\ltimes="2\msb@6E
\mathchardef\rtimes="2\msb@6F \mathchardef\dabar="0\msa@39
\mathchardef\daright="0\msa@4B \mathchardef\daleft="0\msa@4C

\def\Bbb{\ifmmode\let\next\Bbb@\else
        \def\next{\errmessage{Use \string\Bbb\space only in math mode}}\fi\next}
\def\Bbb@#1{{\Bbb@@{#1}}}
\def\Bbb@@#1{\fam\msbfam#1}
\catcode`\@=12



\newcount\senu
\def\senum{\number\senu}
\newcount\ssnu
\def\ssnum{\number\ssnu}
\newcount\fonu
\def\fonum{\number\fonu}

\def\num{{\senum.\ssnum}}
\def\numfo{{\senum.\ssnum.\fonum}}


\outer\def\section#1\par{\vskip0pt
  plus.3\vsize\penalty20\vskip0pt
  plus-.3\vsize\bigskip\vskip\parskip
  \message{#1}\centerline{\sectionfont\senum\enspace#1.}
  \nobreak\smallskip}

\def\endsection{\advance\senu by1\penalty-20\smallskip\ssnu=1}
\outer\def\ssection#1\par{\bigskip
  \message{#1}{\noindent\bf\num\ssectionfont\enspace#1.\thinspace}
  \nobreak\normalfont}

\def\endssection{\advance\ssnu by1\smallskip\ifdim\lastskip<\medskipamount
\removelastskip\penalty55\medskip\fi\fonu=1\normalfont}

\def\proclaim #1\par{\bigskip
  \message{#1}{\noindent\bf\num\enspace#1.\thinspace}
  \nobreak\claimfont}

\def\cor{\proclaim Corollary\par}

\def\lemma{\proclaim Lemma\par}
\def\prop{\proclaim Proposition\par}
\def\rmk{\proclaim Remark\par\normalfont}
\def\thm{\proclaim Theorem\par}

\def\endcor{\endssection}

\def\endlemma{\endssection}
\def\endprop{\endssection}
\def\endrmk{\endssection}
\def\endthm{\endssection}

\def\Proof{{\noindent\sl Proof: \/}}


\def\maplefto#1{\ \smash{\mathop{\longleftarrow}\limits^{#1}}\ }

\def\llongrightarrow{\relbar\joinrel\relbar\joinrel\rightarrow}
\def\lllongrightarrow{\hbox to 40pt{\rightarrowfill}}

\def\twoheadrightarrow{\rightarrow\kern -8pt\rightarrow}

\def\maprighto#1{\smash{\mathop{\longrightarrow}\limits^{#1}}}

\def\llongmaprighto#1{\ \smash{\mathop{\llongrightarrow}\limits^{#1}}\ }

\def\lllongmaprighto#1{\ \smash{\mathop{\lllongrightarrow}\limits^{#1}}\ }

\def\longleftmapsto{\longleftarrow\kern-2pt\mapstochar\;}

\def\llongmapsto{\,\vert\kern-3.2pt\joinrel\longrightarrow\,}
\def\llongmapsto{\,\vert\kern-3.7pt\joinrel\llongrightarrow\,}
\def\lllongmapsto{\,\vert\kern-5.5pt\joinrel\lllongrightarrow\,}

\def\isomarrow{\maprighto{\lower3pt\hbox{$\scriptstyle\sim$}}}
\def\llongisomarrow{\llongmaprighto{\lower3pt\hbox{$\scriptstyle\sim$}}}
\def\lllongisomarrow{\lllongmaprighto{\lower3pt\hbox{$\scriptstyle\sim$}}}

\def\lisomarrow{\maplefto{\lower3pt\hbox{$\scriptstyle\sim$}}}

\font\labprffont=cmtt8
\def\strutdepth{\dp\strutbox}
\def\labtekst#1{\vtop to \strutdepth{\baselineskip\strutdepth\vss\llap{{\labprffont #1}}\null}}
\def\marglabel#1{\strut\vadjust{\kern-\strutdepth\labtekst{#1\ }}}

\def\label #1. #2\par{{\definexref{#1}{\num}{#2}}}
\def\labelf #1\par{{\definexref{#1}{\numfo}{formula}}}
\def\labelse #1\par{{\definexref{#1}{\num}{section}}}


\def\QQ{{\Bbb Q}}

\def\Spec{{\rm Spec}}

\def\SGA1{{\rm SGA1}}

\senu=1 \ssnu=1 \fonu=1

\overfullrule=0pt
\voffset=10mm

 \font\sc=cmcsc9

\font\tenbb=msbm10 \font\sevenbb=msbm7 \font\fivebb=msbm5
\newfam\bbfam
\textfont\bbfam=\tenbb \scriptfont\bbfam=\sevenbb
\scriptscriptfont\bbfam=\fivebb
\def\bb{\fam\bbfam}

\def\O{{\cal O}}
\def\X{{\cal X}}
\def\<{\langle}
\def\>{\rangle}
\def\L{{\cal L}}
\def\I{{\cal I}}
\def\PP{{\Bbb P}}


\centerline{\titlefont Dyson's Theorem for curves.} \spacing
\centerline{\bf  C. Gasbarri}

\bigskip
\bigskip

{\insert\footins{\leftskip\footnotemargin\rightskip\footnotemargin\noindent\eightpoint
$2000$ {\it Mathematics Subject Classification}.  Primary: 11G30,
11G50, 14G40, 14H25.
}

\vbox{{\leftskip\abstractmargin \rightskip\abstractmargin
\eightpoint

\noindent A{{\sixrm BSTRACT}}.\enspace Let $\scriptstyle K$ be a
number field and $\scriptstyle X_1$ and $\scriptstyle X_2$ two
smooth projective curves defined over it. In this paper we prove an
analogue of the Dyson Theorem for the product $\scriptstyle X_1
\times  X_2$. If $\scriptstyle X_i \, = \, {\bb P}_1$ we find the
classical Dyson theorem. In general, it will imply a self contained
and easy proof of Siegel theorem on integral points on hyperbolic
curves and it will give some insight on effectiveness. This proof is
new and avoids the use of Roth and Mordell-Weil theorems, the theory
of Linear Forms in Logarithms and the Schmidt subspace theorem.

}}

\bigskip

\section Introduction\par}

\bigskip

After the proof of the Mordell conjecture by Faltings (the first
proof is in [Fa1], but [Fa2], [B2] and [Vo2] are nearer to the
spirit of this paper), most of the {\it qualitative} results in
the diophantine approximation of algebraic divisors {\it by
rational points} over curves are solved.

Historically, the first concluding result is the Siegel's theorem:
An affine hyperbolic curve contains only finitely many
$S$-integral points; we know that we cannot suppose less on the
{\it geometry} of the involved curve: ${\bb A}^1$ and ${\bb G}_m$
have, as soon as the field is sufficiently big, infinitely many
integral points.

After a long and interesting story of partial results (Liouville,
Thue, Siegel, Dyson, Gelfand$\ldots$), Roth proved that, if
$\alpha$ is an algebraic number then, for every $\kappa > 2$, the
equation
$$
\left\vert \alpha - {p \over q} \right\vert \leq {1 \over \vert q
\vert^{\kappa}}
$$
admits only finitely many solutions ${p \over q} \in {\bb Q}$.
Here again, by Dirichlet's theorem, we know that, for $\kappa =
2$, the equation may have infinitely many solutions.

Eventually, the already quoted theorem of Faltings closes the
story: a compact hyperbolic curve contains only finitely many
rational points.

It is a fact that, from a {\it quantitative} point of view, we are
still very far from a satisfactory answer (up to the very
interesting partial results in [B1], [B3], [BVV] and [BC]): In
each of the three problems quoted above we are not able to give an
upper bound for the heights of the searched solution. And, even
worste, we are not able to say if {\it there is} any solutions to
each of these problems.

Let's have a closer look to the Siegel's theorem: the modern proof
of it relies upon the Roth's theorem and on the Mordell Weil's
theorem or on the theory of the Linear Forms in Logarithms and
again on the Mordell--Weil's theorem; recently, a new proof, based
on the Schmidt's subspace theorem has been given [CZ].
Consequently, if one tries to find an effective proofs by refining
the existing proof, one will crash into the problems of
effectiveness in Roth's theorem and in the computation of a basis
for the Mordell-Weil group of the Jacobian (problem which seems
easier but not yet completely solved) or in the effectiveness in
Schmidt's theorem. Nevertheless some very important cases of
effective Siegel's theorem are given in [Bi]. So, at a first
glance, an effective version of Siegel's theorem will be
consequence of the solutions of other problems, which seems to be
even more difficult. This is very unsatisfactory, also because a
strong effective version of it will imply a version of the
$abc$-conjecture ([Su]).

In this paper we prove a theorem in the spirit of the Dyson's
Theorem [B1] over the product of two curves. It will easily imply
Siegel's theorem. Up to standard facts in algebraic geometry and
in the theory of heights, the theorem is self contained and
essentially elementary. Consequently it release Siegel's theorem
from other big theorems. In this way Siegel's theorem becomes a
result which is completely independent and, perhaps an effective
version of it can be studied on its own.

We now give a qualitative statement of the main theorem of this
paper; for a precise statement, cf. section 2.

Let $K$ be a number field,  let $L_1,\dots L_r$ be  finite
extensions of $K$ and $n := \max\{[L_i\cdot L_j : K]$ and denote
by $A$ the $K$--algebra $\oplus L_i$. Let $X_1$ and $X_2$ be
smooth projective curves over $K$ and $D_i = {\rm Spec} (A)
\rightarrow X_i$, be effective geometrically reduced divisors on
$X_i$; note that the $D_i$'s may have different degrees. Let $H_i$
be a line bundle of degree one over $X_i$ and $h_{H_i} (\cdot)$
height functions associated to $H_i$. Finally, let $S$ be a finite
set of places of $K$ and for every $v\in S$ let $\lambda_{D_i , v}
(\cdot)$ be Weil functions associated to $D_i$ and $v$.

\label teorema1. theorem\par\thm  Let $\vartheta_1$, $\vartheta_2$
and $\epsilon$ be three rational numbers such that $\vartheta_1
\cdot \vartheta_2 \geq 2n+\epsilon$ and $\vartheta_i\geq 1$. Let
$\varphi: S\to [0,1]$ be a function such that $\sum_{v\in
S}\varphi(v)=1$. Then the set of rational points $(P,Q) \in X_1(K)
\times X_2(K)$ such that for every $v\in S$
$$
\lambda_{D_1 , v} (P) > \varphi(v)\cdot\vartheta_1 \cdot h_{H_1} (P)
$$
and
$$
\lambda_{D_2 , v} (Q) >\varphi(v)\cdot \vartheta_2 \cdot h_{H_2} (Q)
$$
is contained in a proper closed subset whose irreducible components
are either fibers or points.
\endthm

If we apply the theorem to ${\bb P}_1 \times {\bb P}_1$ and
$\vartheta_1 = \vartheta_2 = \sqrt{2n} + \epsilon$ we reobtain the
classical theorem of Dyson (cf. [B1]):

\label dyson. corollary\par\cor Let $\alpha$ be an algebraic
number of degree $n$ over ${\bb Q}$. Then there are only finitely
many ${p \over q} \in {\bb Q}$ such that
$$
\left\vert \alpha - {p \over q} \right\vert \leq {1 \over
q^{\sqrt{2n} + \epsilon}} \, .
$$
\endcor

If we apply the theorem to $C \times C$ where $C$ is an arbitrary
curve, $D$ a reduced divisor on it, we obtain the following
generalization:

\label dysonforcurves. corollary\par\cor Let $C$ be a smooth
projective curve over a number field $K$ and $M$ be a line bundle
of degree one on it; let $D$ be a reduced divisor of degree $n$
over $C$ then for all $p \in C(K)$ we have
$$
\lambda_{D,S} (p) \leq ( \sqrt{2n} + \epsilon) \, h_M (p) + O(1)
\, .
$$
\endcor

The involved constant is not effective.

\medskip

Corollary \ref{dysonforcurves} easily implies Siegel's theorem on
$S$--integral points. We first recall the definition of integral
points: let $C$ be a smooth projective curve defined over a number
field $K$. Let $D$ be an effective reduced divisor on $C$. Suppose
that we fixed a logarithmic height function $h_{D}(\cdot)$ with
respect to $D$. Let $S$ be a finite set of places of $K$ and
$\lambda_{D,S}(\cdot)$ be a Weil function associated to $S$ and
$D$ (cf. \S 2 for definitions and references). Let $I\subset C(K)$
be a set of rational points. The set $I$ is said to be {\it
integral with respect to $D$ and $S$} (or $(D,S)$--integral) if
there exists a constant $C$ such that, for every point $P\in I$
$$\vert h_{D}(P)-\lambda_{D,S}(P)\vert\leq C$$
(for short, we will write $\lambda_{D,S}(P)= h_{D}(P)+O(1)$).

\label siegeltheorem. corollary\par\cor (Siegel Theorem) Let $K$ be
a number field and $S$ be a finite set of places of it. Let $C$ be a
smooth projective curve of genus $g$ defined over a number field
$K$. Let $D$ be a reduced effective divisor on $C$ different from
zero. Suppose that $2g-2+\deg(D)>0$. Then every set of
$(D,S)$--integral points is finite.
\endcor

\Proof Fix a line bundle $M$ of degree one on $C$. For every
positive number $\epsilon$, standard properties of heights (cf. for
instance [HS]) give the existence of a constant $A$ such that $\deg
(D)h_M(\cdot)\leq (1+\epsilon)h_D(\cdot)$+A. Suppose that
$\deg(D)\geq 3$. In this case $2g-2+\deg(D)>0$ independently on the
genus. Let $I$ be a set of $(D,S)$--integral points. By definition
$h_{D}(P)=\lambda_{S,D}(P) +O(1)$. Fix $\epsilon_1$ very small and
apply \ref{dysonforcurves}; we  obtain, for every $P\in I$,
$$\eqalign{\deg(D)h_M(P)\leq(1+\epsilon)h_D(P)&=(1+\epsilon)\lambda_{S,D}(P)+0(1)\cr&\leq
(1+\epsilon)(\sqrt{2\deg(D)}+\epsilon_1)h_M(P)+O(1).\cr}$$ If
$\epsilon$ and $\epsilon_1$ are sufficiently small, we have that
$\deg(D)-(1+\epsilon)(\sqrt{2\deg(D)}+\epsilon_1)\geq 0$;
consequently the height, with respect to $M$, of points $P$ in $I$
is bounded independently on $P$. From this we conclude in this
case.

Suppose that $D$ is arbitrary. In this case $g\geq 1$. Take an
\'etale covering $f\colon C'\to C$ of degree bigger then three.
Then $\deg(f^\ast(D))\geq 3$. By the theorem of Chevalley and Weil
([Se] Theorem 4.2) there is a finite extension $K'$ of $K$ such
that $f^{-1}(C(K))\subset C'(K')$. Apply the previous case to
$C'$, $f^\ast(D)$ and $I':=f^{-1}(I)$ and conclude.

\smallskip

Using Roth theorem and the weak Mordell-Weil theorem one obtains,
if $g>0$,
$$
\lambda_{D,S} (p) \leq \epsilon h_M (p) + O(1) \, ;
$$
which is much stronger then \ref{dysonforcurves} (but, it implies
the same qualitative result on integral points). Nevertheless, as
already said, the proof we propose here is much simpler and its
ineffectiveness is essentially self contained: it does not depend
on other theorems.

\smallskip

A remark on the language and the methods used: In this paper we
decided to use the language of arithmetic geometry \`a la
Grothendieck and the Arakelov geometry; although this needs a
little bit of background, which nowadays is (or should be)
standard, this language allows to better understand and compute
the involved constants  and to understand their nature. It is our
opinion that, algebro geometric and Arakelov methods, being more
intrinsic and conceptual, are more adapted to understand the
strategy and the ideas of a proof in diophantine geometry.  In any
case, in the paper we tried to recall the background in Arakelov
geometry needed to understand it. For an introduction to the
Arakelov geometry used in this paper cf. [MB] or the more general
[BGS]. A very fast introduction to the Arakelov geometry of
arithmetic surfaces is in [Ga].

\smallskip

Before we start the proof we would like to give a very informal
argument that roughly explain the strategy of the proof.

Suppose that $D= 0\subset\PP^1$ is the divisor of a single point.
Fix local coordinates $(z_1,z_2)$ around $(0,0)$ in
$\PP^1\times\PP^1$. Let $\pi:\tilde X\to\PP^1\times\PP^1$ be the
blow up of $(0,0)$ and let $E$ be the exceptional divisor.

Suppose that we have a couple of points
$(p_1,p_2)\subset\PP^1\times\PP^1(K)$ such that
$\lambda_{0,\infty}(p_i)\gg h_{\O(1)}(p_i)$(the involved constants
are not important in this argument). In particular we may suppose
that, in the Euclidean topology, $p_i$ is very near to $0$. We
want to prove there are only finitely many couples of such points.

Observe that $\lambda_{0,\infty}(p_i)=-\log{{\vert
z_i(p_i)\vert}\over{\sqrt{1+\vert z_i(p_i)\vert^2}}}\sim
-\log\vert z_i(p_i)\vert$.

Since the exceptional divisor $E$ is locally defined by $z_1$ (or
$z_2$), a local computation gives $h_{\O(E)}(p_1,p_2)\gg
h_{\O(1,1)}(p_1,p_2)$.

In the sequel we denote by $h_i$ the real number $h_{\O(1)}(p_i)$.

Suppose that we can find constants $A_i$ such that for an infinite
sequence of positive integers $d$'s there is a section $f\in
H^0(\PP^1\times\PP^1, \O({{d}\over {h_1}},{{d}\over{h_2}}))$ with
integral coefficients (observe that
$H^0(\PP^1\times\PP^1,\O(d_1,d_2))$ is the space of polynomials
with bidegree $(d_i,d_2)$) such that

-- $\sup\{\Vert f\Vert_{FS}(z_1,z_2)\}\leq
A_1^{{d}\over{h_i+h_2}}$.

-- $div(f)$ vanishes with order at least
$m:={{dA_2}\over{h_1+h_2}}$ on $(0,0)$.

-- $f$ {\it do not vanishes in $(p_1,p_2)$}.

Then the strict transform of $div(f)$ give rise to a section
$\tilde{f}\in H^0(\tilde{X}, \pi^\ast(\O({{d}\over
{h_1}},{{d}\over{h_2}}))-mE)$.

Since $\tilde{f}(p_1,p_2)\neq 0$, we find

$$h_{\pi^\ast(\O({{d}\over{h_1}},{{d}\over{h_2}}))(-mE)}(p_1,p_2)\geq
-{{d}\over{h_1+h_2}}\log(A_1).$$

Thus, since $h_{\O(E)}(p_1,p_2)\gg h_1+h_2$,
$${{d}\over{h_1}}\cdot h_1+{{d}\over{h_2}}\cdot
h_2-{{dA_2}\over{h_1+h_2}}\cdot(h_1+h_2)\geq-{{d}\over{h_1+h_2}}\log(A_1).$$
And from this we deduce that
$$2-A_2\geq -{{\log(A_1)}\over{h_1+h_2}};$$
consequently $p_1$ and $p_2$ must have bounded height.

The general proof need to construct such a section and prove that
it do not vanish on the point. In general one cannot work with the
ideal $(z_1,z_2)$ (in our example we blow up this ideal) but one
consider a more complicated ideal (introduced in \S 3) which
depend on the constants involved in the inequality between Weil
functions and heights we are assuming. To construct the section
with small norm one uses the Siegel Lemma and it will not exist if
we assume a too strong inequality; the inequalities supposed in
theorem \ref{teorema1} allow to construct such a section. One
cannot prove that the section do not vanish on the point, thus one
prove that, under suitable conditions on the heights, the section
has a small order of vanishing on it (this is the more geometrical
part of the proof): this is done in \S 4. This is also the part of
the proof which is not effective. Thus one take a suitable
"derivative of the section" to produce a section which do not
vanish on the point and then the conclusion is essentially the one
explained above.

One should notice that almost all the proofs in diophantine
approximation follow this strategy (for instance, in one take only
one factor, one obtain the Liouville inequality).

\smallskip

I would like to thank the anonymous Referee for her/his comments
and remarks; I could improve and in some cases simplify the
arguments following her/his suggestions and remarks.

\endsection

\bigskip

\section  Statement of the main theorem and notations\par

\bigskip

Let $K$ be a number field and $O_K$ be its ring of integers. We
will denote by $M_K$ the set of (finite and infinite) places of
$K$. Let $M_{\infty}$ be the set of infinite places of $K$. Let
$S$ be a finite subset of $M_K$. We will denote by $O_S$ the ring
of $S$-integers of $K$. For every $v \in M_K$ let $K_v$ be the
completion of $K$ at the place $v$, $O_v$ be the local ring of $v$
and $k_v$ be its residue field. For every scheme $X \rightarrow
{\rm Spec} (O_K)$ we will denote by $X_v$  (resp. $x_{O_v}$, resp.
$X_K$),  the base change of it from ${\rm Spec} (O_K)$ to ${\rm
Spec} (K_v)$ (resp to $\Spec(O_v)$, resp to ${\rm Spec} (K)$).
Similarly, if $L$ is an extension of $K$, we will denote by $O_L$
the ring of integers of $L$, by $X_L$ the base change of $X$ to
$\Spec(L)$ etc.

Let $L_1, \ldots ,L_r$ be finite extensions of $K$ and $O_{L_i}$ be
the ring of integers of $L_i$. We will denote by $A$ the
$O_K$-algebra $\oplus \, O_{L_i}$.

We will denote by $\overline K$ the algebraic closure of $K$.

Let $X\to\Spec(O_K)$ be an $O_K$--scheme. An hermitian vector bundle
$\overline E$ of rank $r$ over $X$ is a couple
$(E,\<\cdot,\cdot\>_{\sigma})_{\sigma\in M_{\infty}}$ where

\noindent -- $E$ is a vector bundle of rank $r$ over $X$.

\noindent -- for every infinite place $\sigma$, the vector bundle
$E_{\sigma}$  is an holomorphic vector bundle over the ${\bb
C}$--scheme $X_{\sigma}$; then $\<\cdot,\cdot\>_{\sigma}$ is a
$C^{\infty}$ metric on $E_{\sigma}$ (and if
$\tau=\overline{\sigma}$, the metric on $E_\tau$ is the complex
conjugate of the metric on $E_\sigma$).

If $M$ is an hermitian vector bundle of rank one, we will call it
{\it hermitian line bundle}. If $M$ is an hermitian line bundle over
$\Spec(O_K)$ we will define its Arakelov degree by the following
formula: Let $s\in M\setminus\{ 0\}$; then
$$\widehat{\deg}(M):=\log({\rm Card}(M/s\cdot O_K))-\sum_{\sigma\in
M_{\infty}}\log\Vert s\Vert_\sigma.$$ This formula is well defined
because of the product formula (cf. for instance [SZ]).

If $\overline E$ is an arbitrary hermitian vector bundle over
$\Spec (O_K)$ then the line bundle $\bigwedge^{max}E$ is
canonically equipped with an hermitian metric; consequently we can
define the hermitian line bundle $\bigwedge^{\max}\overline E$. We
then define $\widehat{\deg}(\overline
E):=\deg(\bigwedge^{\max}(\overline E))$.

Suppose that $\overline E_1$ and $\overline E_2$ are hermitian
vector bundles over $\Spec(O_K)$ and $f\colon E_1\to E_2$ is a
linear map. Then, for every infinite place $\sigma$, $f$ induces a
linear map $f_\sigma\colon (E_1)_{\sigma}\to (E_2)_{\sigma}$; Let
$\Vert f_\sigma\Vert_\sigma$ be the norm of it. Then we define
$\Vert f\Vert:=\sup\{\Vert f_{\sigma}\Vert_{\sigma}\}_{\sigma\in
M_{\infty}}$.

More generally: Suppose that $X\to\Spec(O_K)$ is an arithmetic
scheme and $\overline E$ is an hermitian vector bundle over it.
Suppose that, for every $\sigma\in M_{\infty}$ the complex variety
$X_{\sigma}(\Bbb C)$ is projective and smooth and that we fixed a
smooth hermitian metric on it. Under these conditions the
$O_K$--module $H^0(X,E)$ has a natural structure of hermitian
$O_K$--module: indeed, for every $\sigma\in M_{\infty}$, the
complex vector space $H^0(X,E)_{\sigma}$ is equipped with the
$L^2$ hermitian metric induced by the metric on $X_{\sigma}(\Bbb
C)$ and on $E_{\sigma}$. For every infinite place $\sigma$,  the
complex vector space $H^0(X,E)_{\sigma}$ is naturally equipped
with the $\sup$ norm: $\Vert f\Vert_{\sup, \sigma}:=\sup_{x\in
X_{\sigma}(\Bbb C)}\{ \Vert f\Vert(x)\}$. The $L^2$ and $\sup$
norms are comparable (as explained for instance in [Bo]);
consequently we can work with the norm we prefer.

Let $f_1 : {\cal X}_1 \rightarrow B := {\rm Spec} (O_K)$ and $f_2
: {\cal X}_2 \rightarrow B$ be two regular, semistable arithmetic
surfaces over $O_K$. Let $\Delta_i \hookrightarrow {\cal X}_i
\times_B {\cal X}_i$ ($i = 1,2$) be the diagonal divisor. The
divisor $\Delta_i$ is, a priori, just a Weil divisor (the scheme
${\cal X}_i\times{\cal X}_i$ may be not regular); let $\tilde
{{\cal X}_i\times{\cal X}_i}$ be the blow up of it along
$\Delta_i$ and $\tilde{\Delta_i}$ be the exceptional divisor.

For every infinite place $\sigma$, we fix a symmetric hermitian
structure on the line bundle $({\cal O}
(\tilde{\Delta_i}))_{\sigma}$ ($i = 1,2$). Let $\sigma \in M_K$ be
an infinite place and $P \in ({\cal X}_i)_{\sigma} ({\bb C})$;
denoting by $\iota_P : ({\cal X}_i)_{\sigma} ({\bb C}) \rightarrow
({\cal X}_i \times {\cal X}_i)_{\sigma} ({\bb C})$ the embedding
deduced from the map $\iota_P (x) := (x,P)$, we have a canonical
isomorphism $\iota_P^* {\cal O} (\Delta) \simeq {\cal O} (P)$. For
every place $\sigma$ and $P\in({\cal X}_i)_{\sigma} ({\bb C})$, we
put on $\O(P)$ the metric obtained taking the pull back metric via
$\iota_p$. As a consequence, for every divisor $D$ of ${\cal
X}_i$, the line bundle ${\cal O} (D)$ is equipped with a canonical
metric (depending only on the choices made until now).

Let $D$ be an effective divisor on $(\X_i)_K$. For every finite set
of places $S\in M_K$ we can choose a canonical representative for
the Weil function $\lambda_{D,S}(\cdot)$ in the following way: First
of all we take the schematic closure of $D$ on $\X_i$; this will be
a Cartier divisor over $\X_i$.

\noindent -- Suppose that $S:=\sigma$ is an infinite place; let
${\bb I}_D$ be the canonical section of $(\O(D))_{\sigma}$. Let
$\Vert\cdot\Vert(\cdot)$ be the metric on ${\O(D)}_\sigma$ defined
above; then we define, for every $x\in (\X_i)_{\sigma}(\Bbb
C)\setminus\{D\}$:
$$\lambda_{D,\sigma}(x):=-\log\Vert{\bb I}_D\Vert(x)$$

\noindent -- Suppose that $S:=v$ is a finite place. Since $D$ and
$(\X_i)_v$ are generic fibers of their models over $\Spec(O_v)$, as
explained in [D], the line bundle $(\O(D))_{v}$ over the
$K_v$--scheme $(\X_i)_v$ is equipped with a $v$--adic norm;
consequently the Weil function $\lambda_{D,v}(\cdot)$ is defined
similarly.

\noindent -- If $S$ is arbitrary, then $\lambda_{D,S}(\cdot)$ is
defined as sum if local terms as explained for instance in [HS]
chapter B 8.

The choice of a metric on the diagonal induces a metric on the
relative dualizing sheaf $\omega_{{\cal X}_i / B}$; we fix such a
metric; remark that, by construction, the adjunction formula
holds: for every section $P:B \rightarrow {\cal X}_i$ we have a
canonical isomorphism \labelf adjunctionformula\par $$
\omega_{{\cal X}_i / B} \vert_P \simeq {\cal O} (-P) \vert_P
\eqno{{(\numfo)}}$$\advance\ssnu by1 of {\it hermitian} line
bundles on B. For a general reference on this cf. [MB]. For a
reference on Weil functions cf. [HS].

For every hermitian line bundle $\overline M:=(M;\Vert\cdot\Vert)$
over $\X_i$ we can define a height function
$$h_{\overline M}(\cdot)\colon (\X_i)_K(\overline K)\longrightarrow{\bb R}$$
in the following way:

\noindent Let $P\in \X_i(\overline K)$. It is defined over a
finite extension $L$. Let $(\X_i)_{O_L}\to\Spec(O_L)$ be the
minimal regular model of $(\X_i)_L$. The point $P$ corresponds to
a section ${\cal P}\colon \Spec(O_L)\to(\X_i)_{O_L}$; we define
$$h_M(P):={{1}\over{[L:{\bb Q}]}}\cdot\deg({\cal P}^\ast(\overline
M)).$$

An hermitian line bundle $\overline M$ on $\X_i$ is said to {\it
be arithmetically ample} if its degree on the projective curve
$X_K$ is positive and $h_M(\cdot)>0$.

Fix an arithmetically ample hermitian line bundles $(M_i, \Vert
\cdot \Vert_{M_i})$ on ${\cal X}_i$ of generic degree one.

We will denote by $(\cdot ; \cdot)$ the Arakelov intersection
pairing on each of the ${\cal X}_i$ as defined for instance in [BGS]
or [MB].

If $D$ is an effective reduced divisor over ${\cal X}_i$; write $D
:= \sum D_j$ where each $D_j$ is an irreducible divisor. Define
the following three numbers associated to it:

\noindent -- Let $L_j$ be an extension of $K$ where $D_j$ splits
as sum of points: if $f_j:(\X_i)_{O_{L_j}}\to \X_j$ is the minimal
regular model of the base change of $\X_i$ to $\Spec(O_{L_j})$,
then $f_j^{\ast}(D_j)=\sum_h P_{hj}+V$; where $P_{hj}$ are
sections and $V$ is a vertical divisor: Then we define
$$S(D):=\max_{h \; ,\; j}\{ -{{1}\over{[L_j:\Bbb Q]}}\cdot
(\O(P_{hj});\O(P_{hj})); 1\};$$
$$H(D):=\max_{h \; ,\; j}\{ h_{M_i}(P_{ij});1\};$$
and
$$T(D):=S(D)\cdot H(D).$$
Observe that, by formula \ref{adjunctionformula}, we have that
$-(\O(P_{hj});\O(P_{hj}))=(\omega_{\X_i/O_L};\O(P_{hj}))$.

We eventually fix a positive integer and three positive rational
numbers $\vartheta_1, \vartheta_2$ and $\epsilon$ such that
$$
\vartheta_1\cdot \vartheta_2 \geq 2n + \epsilon \, .
$$

The main theorem of this paper is the following generalization of
Dyson's theorem:

\label maindyson. theorem\par\thm Under the hypotheses above there
exist two effectively computable constants $R_1$ and $R_2$,
depending only on the ${\cal X}_i$, the hermitian line bundles
$M_i$, the metrics on the diagonals, the $\vartheta_i$ and the
constant $\epsilon$, for which the following holds:

Let $L_1 , \ldots , L_r$ be finite extensions of $K$; denote by $n$
the number $n := {\rm max}\{ [ L_i \cdot L_j : K]\}$, by $O_{L_i}$
the ring of integers of $L_i$ and by ${\cal A}$ the $O_K$-scheme
${\cal A} := {\rm Spec} (\oplus \, O_{L_i})$. Let $\varphi\colon
S\to [0,1]$ be a function such that $\sum_{v\in S}\varphi (v)=1$.

Let
$$
D_i : {\cal A} \rightarrow {\cal X}_i
$$
be reduced effective divisors over ${\cal X}_i$ ($i = 1, 2$).

If $(P,Q) \in {\cal X}_1(K) \times {\cal X}_2 (K)$ is a couple of
rational points such that

\smallskip

\item{(a)} $h_{M_1} (P) \geq R_1 \cdot T(D_1) \cdot T(D_2)$

\smallskip

\item{(b)} for every $v\in S$
$$\lambda_{D_1,v} (P) >\varphi(v)\cdot \vartheta_1 \cdot h_{M_1} (P)\;\;
{\rm and}\;\; \lambda_{D_2,v} (Q) > \varphi(v)\cdot\vartheta_2 \cdot
h_{M_2} (Q);$$

\smallskip

\noindent then
$$
h_{M_2} (Q) \leq R_2 \cdot T(D_1) \cdot T(D_2) \cdot h_{M_1} (P)
\, .
$$
\endthm

\medskip

This will easily imply the qualitative theorem and its
corollaries.

In the following sections we will introduce the tools we need for
the proof of \ref{maindyson}, we will give it in the last section.

\endsection

\vskip 1cm

\section Small sections\par

\bigskip

Let $L$ be a finite extension of $K$ of degree $n$ and $O_L$ its
ring of integers. Denote by $B_L$ the scheme $\Spec(O_L)$.

Let ${\cal L}$ be a line bundle over $B := {\rm Spec} (O_K)$; we
will denote by ${\cal O}[{\cal L}]$ the $O_K$-algebra ${\rm Sym}
(\oplus \, {\cal L}^{\otimes n})$ and by ${\cal O}[\![{\cal
L}]\!]$ the $O_K$-algebra $\prod {\cal L}^{\otimes n}$ with the
multiplicative structure given by $(a_n) \cdot (b_n) := (c_n)$
where $c_n := \sum_{i+j=n} \,  a_i \otimes b_j$ (if ${\cal L}$ is
the trivial line bundle ${\cal O}_B$ then ${\cal O} [\![{\cal
O}_B]\!]$ is the usual ring of power series $O_K [\![X]\!]$). If
${\cal L}_1$ and ${\cal L}_2$ are two line bundles we define
${\cal O} [{\cal L}_1, {\cal L}_2]$ and ${\cal O} [\![ {\cal L}_1,
{\cal L}_2]\!]$ in a similar way.

Let $f_{\cal L} : {\bb V} ({\cal L}) \rightarrow B$ be the affine
$B$-scheme ${\rm Spec} ({\cal O} [{\cal L}])$ then it is easy to
verify that:

\smallskip

\item{(a)} there is a canonical isomorphism $f^* ({\cal L}) \simeq \Omega_{{\bb V}({\cal L}) / B}^1$;

\smallskip

\item{(b)} if ${\bf 0} : B \rightarrow {\bb V} ({\cal L})$ is the canonical section,
there is a canonical isomorphism
$\widehat{{\bb V} ({\cal L})}_{\bf 0} \simeq {\rm Spf} ({\cal O} [\![ {\cal L} ]\!])$.

\smallskip

Suppose that $\overline\L_1$ and $\overline\L_2$ are hermitian
line bundles. Let $\sigma\in M_{\infty}$. For every positive
integer $n$, the complex vector space
$\bigoplus_{a+b=n}(\L_1^{\otimes a}\otimes\L_2^{\otimes
b})_\sigma$ has a natural structure of hermitian vector space.
Consequently also $(\O[\L_1,\L_2])_\sigma=\bigoplus_{n\geq
0}\bigoplus_{a+b=n}(\L_1^{\otimes a}\otimes\L_2^{\otimes
b})_\sigma$ has a natural structure of hermitian vector space. Let
$J\subset\O[\L_1,\L_2]$ be an ideal; since
$(\O[\L_1,\L_2])_\sigma$ is direct sum of finite dimensional
hermitian vector space, we can find an orthonormal basis ${\cal
B}_{\sigma}$ of $(\O[\L_1,\L_2])_\sigma$ such that ${\cal
B}_{\sigma}$ is disjoint union of ${\cal B}_1$ and ${\cal B}_2$
with ${\cal B}_1$ orthonormal basis of $J_\sigma$. Consequently
the vector space $(\O[\L_1,\L_2]/J)_\sigma$ is canonically (via
the projection) isomorphic to $J_{\sigma}^\perp$, thus it is
equipped with the structure of hermitian vector space. Moreover,
suppose that $J_1\subset J_2$, then the metric induced by the
canonical projection $\O[\L_1,\L_2]/J_2\to \O[\L_1,\L_2]/J_1$ is
the given metric.

Let $f : {\cal X} \rightarrow {\rm Spec} (O_K)$ be an arithmetic
surface as in the previous section and let $D : {\rm Spec} (O_L)
\rightarrow {\cal X}$ be a reduced divisor over ${\cal X}$.

Let $f_L : {\cal X}_L \rightarrow {\rm Spec} (O_L)$ be a
desingularization of the arithmetic surface ${\cal X} \times_B
{\rm Spec} (O_L)$. The base change of the morphism $D$ give rise
to a section $S_D : B_L \rightarrow {\cal X}_L$; moreover, if $p :
{\cal X}_L \rightarrow {\cal X}$ is the natural projection, by
construction we have that $p \circ S_D = D$.

\label completiononasection. proposition \par\prop Let
$(\widehat{{\cal X}_L})_D$ be the completion of ${\cal X}_L$ around
$S_D(B_L)$; then there is a natural isomorphism
$$
\Psi_D : (\widehat{{\cal X}_L})_D \longrightarrow {\rm Spf} ({\cal
O} [\![ {\cal O} (-S_D) \vert_{S_D} ]\!]) \, .
$$
\endprop

\Proof Since $\X_L$ is regular and $S_D$ is a section,  $S_D
(B_L)$ is contained in the smooth open set of the structural
morphism $f_L$. Consequently we can find an open neighborhood
${\cal U}$ of $S_D (B_L)$ in ${\cal X}_L$ and an {\it \'etale} map
$g_D : {\cal U} \rightarrow {\bb V} ({\cal O} (-S_D)) \vert_{S_D}$
sending $S_D (B_L)$ to the zero section. From this the proposition
follows.

Let $\X_i$ ($i=1,2$) be the arithmetic surfaces fixed in the
previous section. Let $D_1 : {\rm Spec} (O_{L_1}) \rightarrow
{\cal X}_1$ and $D_2 : {\rm Spec} (O_{L_2}) \rightarrow {\cal
X}_2$ be effective reduced divisors on ${\cal X}_1$ and ${\cal
X}_2$ respectively; let $L := L_1 \cdot L_2$ be the composite of
$L_1$ and $L_2$ over $K$. As before they define two sections $S_i
: {\rm Spec} (O_L) \rightarrow ({\cal X}_i)_{O_L}$ ($i = 1, 2$), .
Let $\xi_{D_1 , D_2} : B_L \rightarrow ({\cal X}_1 \times {\cal
X}_2)_L$ be the point obtained from $S_1$ and $S_2$. and denote by
$(\widehat{{\cal X}_1 \times {\cal X}_2})_{\xi_{D_1 , D_2}}$ the
completion of $({\cal X}_1 \times {\cal X}_2)_L$ around $\xi_{D_1
, D_2}$. As corollary of \ref{completiononasection} we obtain:

\label completionontwopoints. corollary\par\cor Let $(\widehat{{\cal
X}_1 \times {\cal X}_2})_{\xi_{D_1 , D_2}}$ the completion of
$({\cal X}_1 \times {\cal X}_2)_L$ around $\xi_{D_1 , D_2}$. Then
there is a natural isomorphism
$$
\Psi_{D_1 , D_2} : (\widehat{{\cal X}_1 \times {\cal
X}_2})_{\xi_{D_1 , D_2}} \longrightarrow {\rm Spf} ({\cal O} [\![
{\cal O} (-S_1) \vert_{S_1} ; {\cal O} (-S_2) \vert_{S_2} ]\!]) \, .
$$
\endcor

Let $M_L$ be the set of places of $L$; and $\sigma \in M_L$ be an
infinite place. As explained before, the $O_L$-algebra $({\cal O}
[({\cal O} (-S_1) \vert_{S_1} ; {\cal O} (-S_2)
\vert_{S_2}])_{\sigma}$ is naturally equipped with the structure of
{\it hermitian} vector space because of the choice of the metrics as
in the first section.

If $p_i : ({\cal X}_1)_L \times ({\cal X}_2)_L \rightarrow ({\cal
X}_i)_L$ is the natural projection, and $N$ is a line bundle on
$({\cal X}_i)_L$, by abuse of notation, we will denote again by $N$
the line bundle $p_i^* (N)$ on $({\cal X}_1)_L \times ({\cal
X}_2)_L$.

In this section we will construct sections of small norm of
suitable line bundles with high order of vanishing along
$\xi_{1,2}$. As usual the key lemma is the Siegel Lemma. Before we
give the statement (and the proof) of the Siegel Lemma we need, we
recall without proof all the tools we need; for the proofs we
refer to [Bo] \S4.1 and [Sz]:

\smallskip

\item{a)} If $E$ is an hermitian vector bundle over $O_K$, then we call the real number
$\mu_n (E) := {1 \over [K ; {\bb Q}]} \cdot {\widehat{\rm deg} (E) \over rk(E)}$, {\it the slope} of $E$;

\smallskip

\item{b)} within all the sub bundles of a given hermitian vector bundle $E$, there is one $F$ having
maximal slope; we call its slope {\it the maximal slope} of $E$
and denote it by $\mu_{\rm max} (E)$; if $F=E$ we will say that
$E$ is {\it semistable}; by construction the sub bundle $F$ is
semistable;

\smallskip

\item{c)} if $E_1$ and $E_2$ are two hermitian vector bundles, we have that
$\mu_{\rm max} (E_1 \oplus E_2) = {\rm max} \{ \mu_{\rm max} (E_1)$; $\mu_{\rm max} (E_2) \}$;

\smallskip

\item{d)} let $f : E \rightarrow F$ be an {\it injective} morphism between hermitian vector bundles; then
${{1}\over{[K:\bb Q]}}\widehat{\rm deg} (E) \leq rk (E) (\mu_{\rm
max} (F) + \log \Vert f \Vert)$;

\smallskip

\item{e)} there is a constant $\chi (K)$ depending only on $K$ (for the precise value we refer to [Sz])
such that, if $E$ is an hermitian vector bundle on $K$ with
$\widehat{\rm deg} (E) > -rk(E) \chi(K)$, then there is a non
torsion element $v \in E$ such that, for every infinite place
$\sigma$ we have $\sup_{\sigma\in M_{\infty}}\{\log(\Vert v
\Vert_{\sigma})\} \leq 3\log(rk(E))$; we define $\Vert \cdot
\Vert_{\rm sup}$ to be ${\rm sup} \{ \Vert \cdot \Vert_{\sigma}
\}_{\sigma \in M_{\infty}}$ (cf. [BGS] thm 5.2.4 and below it);

\smallskip

\item{f)} let $M_{\infty}$ be the set of infinite places of $K$ and
$\lambda := (\lambda_{\sigma})_{\sigma \in M_{\infty}}$ be an element
of ${\bb R}^{[K:{\bb Q}]}$ with $\lambda_{\sigma} = \lambda_{\overline{\sigma}}$;
we denote by ${\cal O} (\lambda)$ the hermitian line bundle
$(O_K , \Vert 1 \Vert_{\sigma} = \exp (-\lambda_{\sigma}))$.
If $E$ is an hermitian vector bundle over $O_K$ then we denote
by $E(\lambda)$ the hermitian vector bundle $E \otimes {\cal O} (\lambda)$.

\smallskip

\item{g)} (Hilbert-Samuel Formula) there is a constant $C$, depending on the choices made
(but not on the $d_i$'s), such that, if $d_1$ and $d_2$ are
sufficiently big, the Hermitian $O_K$-module $H^0 = ({\cal X}_1
\times {\cal X}_2 , M_1^{d_1} \otimes M_2^{d_2})$ is generated by
elements of sup-norm, less or equal then $C^{d_1 + d_2}$.

\smallskip

We will also need the following

\label exactsequence. lemma\par\lemma Let
$$
0 \rightarrow E_1 \longrightarrow E \longrightarrow E_2
\rightarrow 0
$$
be an exact sequence of hermitian vector bundles; then
$$
\mu_{\rm max} (E) \leq {\rm max} \{ \mu_{\rm max} (E_1) , \mu_{\rm
max} (E_2) \} \, .
$$
\endlemma
\medskip

The proof is straightforward and left to the reader.

Let $K\subseteq L$ be a finite extension and
$\Spec(O_L)\to\Spec(O_K)$ the induced morphism.  Let $F$ be an
hermitian vector bundle on $\Spec(O_L)$. Observe that the vector
bundle $f_\ast(F)$ over $\Spec(O_K)$ is naturally equipped with
the structure of hermitian vector bundle.

\label filtration1. lemma\par\lemma Suppose that $F$ is equipped
with a filtration $F=F_0\supseteq F_1\supseteq
F_2\supseteq\cdots\supseteq F_N=0$ with $F_i/F_{i+1}$ line
bundles, equipped with the induced hermitian metric. Then
$$\mu_{\max}(f_\ast(F))\leq {{1}\over{[L:\QQ]}}\max\{\widehat{\deg}(F_i/F_{i+1})\}.$$
\endlemma

\Proof  By devissage we are reduced to prove it when $F$ is itself
a line bundle. Let $Q\subseteq f_\ast(F)$ be the maximal
semistable subbundle. We deduce a map $f^\ast(Q)\to F$
consequently an isometric inclusion of $O_L$ in
$f^{\ast}(Q^\vee)\otimes F$. Thus we get
$\mu_n(f^\ast(Q^\vee)\otimes F)\geq 0$ {\it because
$f^{\ast}(Q^\vee)$ is semistable}. So $[L:\QQ]\mu_n(Q)\leq
\deg(F)$. the conclusion follows.

\smallskip

The Siegel Lemma we need is the following:

\bigskip

\label siegellemma. lemma\par\lemma (Siegel Lemma) Let $V$ and $W$
be hermitian vector bundles over ${\cal O}_K$. Let $\gamma : V
\rightarrow W$ be a {\rm non injective} morphism. Let $m = rk(V )$
and $n := rk({\rm Ker} (\gamma))$. Suppose that there is a constant
$C > 1$ such that:

\smallskip

\item{i)} $V$ is generated by elements with $\sup$ norm at most $C$;

\smallskip

\item{ii)} $\Vert\gamma\Vert\leq C$;

\smallskip

\noindent then, there is a non zero element $x \in {\rm Ker}
(\gamma)$ with
$$
\sup_{\sigma \in M_{\infty}} \{ \log (\Vert x \Vert_{\sigma} ) \}
\leq {m \over n} \cdot \log (C^2) + \left( {m \over n} - 1 \right)
\mu_{\rm max} (W) - {\chi (K) \over [K:{\bb Q}]}+3\log(n) \, .
$$
\endlemma
\medskip

\Proof Denote by $U$ the hermitian vector bundle ${\rm Ker}
(\gamma)$ with the induced metric. Observe that, by property (e)
above, if $\widehat{\rm deg} (U(\lambda)) > -n\chi(K)$, then there
is a non torsion element $x \in U$ such that
$$
\sup_{\sigma \in M_{\infty}} \{ \log ( \Vert x \Vert_{\sigma}) \}
\leq \sup_{\sigma \in M_{\infty}} \{ \lambda_{\sigma}
\}+3\log(rk(U)) \, .
$$

An easy computation gives $\widehat{\rm deg} (U(\lambda)) =
\widehat{\rm deg} (U) + n \cdot \sum_{\sigma} \lambda_{\sigma}$. Let
$W'$ be the image of $\gamma$. Put on $W'$ the metric induced by the
surjection. Thus we have
$$
\widehat{\rm deg} (U(\lambda)) = \widehat{\rm deg} (V) -
\widehat{\rm deg} (W') + n \cdot \sum_{\sigma} \lambda_{\sigma} \,
.
$$
By property (d) we have ${\widehat{\rm deg} (W') \over [K:{\bb Q}]}
\leq (m-n) (\mu_{\rm max} (W) + \log (C))$ and by the very
definition of Arakelov degree, $\widehat{\rm deg}(V ) \geq -m
[K:{\bb Q}] \log (C)$. Consequently
$$
\eqalign{ \widehat{\rm deg} (U(\lambda)) = & \ \widehat{\rm deg}
(V) - \widehat{\rm deg} (W') + n \cdot \sum_{\sigma \in
M_{\infty}} \lambda_{\sigma} \cr \geq &-2m [K : {\bb Q}] \log (C)
- (m-n) [K : {\bb Q}] \mu_{\rm max} (W) + n \cdot \sum_{\sigma \in
M_{\infty}} \lambda_{\sigma} \, ; \cr }
$$
thus, take $\lambda_{\sigma} = {m \over n} \cdot \log (C^2) +
\left( {m \over n} -1 \right) \mu_{\rm max} (W) - {\chi (K) \over
[K:{\bb Q}]} + \epsilon$ and apply the observation above. The
conclusion follows.

\smallskip

Let $\vartheta_1$, $\vartheta_2$ and $\delta$ be three positive
rational numbers. For every couple of positive integers $(d_1, d_2)$
we denote by ${\cal I}_{\underline{\vartheta} , \delta ,
\underline{d}}$ the ideal sheaf of $({\cal X}_1)_L \times ({\cal
X}_2)_L$ defined by \labelf ideal\par$$ \sum_{{i \over d_1} \cdot
\vartheta_1 + {j \over d_2} \cdot \vartheta_2 \geq \delta \atop i
\leq d_1 , j \leq d_2} {\cal O} (-i S_1) \otimes {\cal O} (-j S_2)
\, . \eqno{{(\numfo)}}
$$
In the same way, we will denote by $I_{\underline{\vartheta} ,
\delta , \underline{d}}$ the ideal of ${\cal O} [{\cal O} (-S_1)
\vert_{S_1} , {\cal O} (-S_2) \vert_{S_2}]$ defined by a condition
analogous to condition \ref{ideal}.

We denote by ${\cal A}_{\underline{\vartheta} , \delta ,
\underline{d}}$ the subscheme of $({\cal X}_1)_L \times ({\cal
X}_2)_L$ defined by the ideal ${\cal I}_{\underline{\vartheta} ,
\delta , \underline{d}}$ and by $W_{\underline{\vartheta} , \delta ,
\underline{d}}$ the $O_L$ algebra ${\cal O} [{\cal O} (-S_1)
\vert_{S_1} , {\cal O} (-S_2) \vert_{S_2}] /
I_{\underline{\vartheta} , \delta , \underline{d}}$. Then:

\smallskip

\item{(i)} the isomorphism $\Psi_{D_1^{h_1} , D_2^{h_2}}$ induces an isomorphism
$\Psi_{1,2} : {\cal A}_{\underline{\vartheta} , \delta , \underline{d}} \rightarrow {\rm Spec} (W_{\underline{\vartheta} , \delta , \underline{d}})$;

\smallskip

\item{(ii)} The $O_L$ module $W_{\underline{\vartheta} , \delta , \underline{d}}$ has
a natural structure of hermitian $O_L$-module. Moreover the
$O_L$-module $W_{\underline{\vartheta} , \delta , \underline{d}}$
has a filtration by $O_L$-submodules $F^{\mu}$ such that $F^{\mu} /
F^{\mu + 1} \simeq {\cal O} (-i S_1) \vert S_1 \otimes {\cal O} (-j
S_2) \vert_{S_2}$ with ${i \over d_1} \cdot \vartheta_1 + {j \over
d_2} \cdot \vartheta_2 \leq \delta$; this filtration is isometric.

\bigskip

\label smallsection. proposition\par\prop Let $\epsilon
> 0$  and $\delta>2$ be  two rational numbers; suppose that
$2\cdot\vartheta_1 \cdot \vartheta_2 > \delta^2[L : K] +
\epsilon$, then there exists a constant $A$ depending only on
${\cal X}_i$, $M_i$, $[L:K]$, $\vartheta_i$ and $\epsilon$ such
that the following holds:

For every couple of irreducible divisor $D_1 \hookrightarrow {\cal
X}_1$, and $D_2 \hookrightarrow {\cal X}_2$ as above and every
couple of sufficiently big integers $(d_1, d_2)$, there is a non
zero section $f \in H^0 ({\cal X}_1 \times {\cal X}_2 , M_1^{d_1}
\otimes M_2^{d_2})$ vanishing along ${\cal
A}_{\underline{\vartheta} , \delta , \underline{d}}$ and such
that, for every infinite place $\sigma \in M_K$ we have
$$
\log (\Vert f \Vert_{\sigma}) \leq {A \over \epsilon} \cdot T(D_1)
\cdot T(D_2) \cdot (d_1 + d_2) \, ;
$$
where the $T(D_i)$ are defined as in \S2.\endprop

\Proof Let $\gamma:\Spec(O_L)\to \Spec(O_K)$ the morphism induced
by the inclusion $K\subseteq L$. It induces a morphism of
hermitian modules
$$\gamma_{d_1 , d_2} : H^0 ({\cal X}_1 \times {\cal X}_2 ,
M_1^{d_1} \times M_2^{d_2}) \rightarrow
\gamma_\ast(W_{\underline{\vartheta} , \delta , \underline{d}}
\otimes (M_1^{d_1}) \vert_{S_1} \otimes (M_2^{d_2})
\vert_{S_2}).$$ Let $K(d_1, d_2)$ be the kernel of $\gamma_{d_1 ,
d_2}$. We have to prove that there exists an element in $K(d_1,
d_2)$ having bounded norm.

In the sequel of this proof, "absolute constant" will be
equivalent to say "a constant which depends only on the $\X_i$, on
the hermitian line bundles $\overline M_i$ and on the metrics on
the diagonals; but {\it independent on the $D_i$'s} and on the
$d_i$'s".

By lemma \ref{exactsequence},  lemma \ref{filtration1} and (ii)
above we can find an absolute constant  $A$ such that $\mu_{\rm
max} (\gamma_\ast( W_{\underline{\vartheta} , \delta ,
\underline{d}} \otimes (M_1^{d_1}) \vert_{S_1} \otimes (M_2^{d_2})
\vert_{S_2})) \leq A \cdot T(D_1)\cdot T(D_2)\cdot (d_1+d_2)$.

By (g) above, as soon as $d_1$ and $d_2$ are sufficiently big, we
can find an absolute constant $A$ for which $H^0 ({\cal X}_1
\times {\cal X}_2 , M_1^{d_1} \otimes M_2^{d_2})$ is generated by
elements with norm bounded by $A^{d_1 + d_2}$.

Now we come to the main part of the proof: we can find an absolute
constant $C$ for which the $O_K$-module $H^0 ({\cal X}_1 \times
{\cal X}_2 , M_1^{d_1} \otimes M_2^{d_2})$ has rank which is
bounded below by $C \cdot d_1 \cdot d_2$. The rank of the
$O_L$-module $W_{\underline{\vartheta} , \delta , \underline{d}}
\otimes (M_1^{d_1}) \vert_{S_1} \otimes (M_2^{d_2}) \vert_{S_2}$
can be bounded from above as follows: the number of the terms of
the filtration described in (ii) is the number of couples of
positive integers $(i, j$) with $i \leq d_1$, $j \leq d_2$ and ${i
\over d_1}\cdot\vartheta_1 + {j \over d_2}\cdot\vartheta_2 \leq
\delta$; as soon as $d_1$ and $d_2$ are sufficiently big, this
number is bounded above by $d_1 \cdot d_2$ multiplied by the area
of the triangle with vertices $(0, 0)$,
$({{\delta}\over{\vartheta_1}} , 0)$, and $(0 ,
{{\delta}\over{\vartheta_2}})$ plus a very small error term,
consequently
$$
rk_{O_K} \left(W_{\underline{\vartheta} , \delta , \underline{d}}
\otimes (M_1^{d_1}) \vert_{S_1} \otimes (M_2^{d_2})
\vert_{S_2}\right) \leq d_1 \cdot d_2 \cdot {\delta^2 \over 2
\vartheta_1 \cdot \vartheta_2} [L : K] + \epsilon' \, .
$$
Consequently there is an absolute constant $A$ such that
$$
{rk_{O_K} \left( h^0 ({\cal X}_1 \otimes {\cal X}_2 , M_1^{d_1}
\times M_2^{d_2})\right) \over rk_{O_K} (K (d_1 , d_2))} \leq {A
\over \epsilon} \, .
$$

For every infinite place $\sigma$ of $K$, we cover the Riemann
surface ${\cal X}_{i,\sigma}$ with a finite number of disks over
which the line bundle $M_i$ trivializes; inside each disk we take
a disk with same center and radius one half of the radius of it;
we may suppose that also these smaller disks cover the Riemann
surface (we suppose that this covering is fixed once for all, in
particularly independently of the $D_i$'s). From the lemma
\ref{cauchylemma} below we deduce that we can find a constant $A$,
independent on the $D_i$'s, such that for every infinite place
$\sigma$ we have $\Vert \gamma_{d_1 , d_2} \Vert_{\sigma} \leq
A^{d_1 + d_2}$. We may suppose that that the $d_i$'s are so big
that $\log(d_i)\leq \epsilon d_i$. We apply now \ref{siegellemma}
to this situation and conclude the proof of the proposition.

\label cauchylemma. lemma\par\lemma  Let $\Delta_R$ be the disk of
radius $R$. Let $f(x,y)$ be an holomorphic function on $\Delta_R
\times \Delta_R$ and $(z_1, z_2) \in \Delta_{R/2} \times
\Delta_{R/2}$ then for every $(i, j)$
$$
\left\vert {\partial^{i+j} f \over \partial x^i \partial y^j} (z_1
, z_2) \right\vert \leq {2^{i+j} i! j! \over R^{i+j}} \cdot
\max_{\Vert x \Vert = \Vert y \Vert = R} \{ \vert f(x,y) \vert \}
\, .
$$
\endlemma

The proof of the Lemma is a straightforward application of the
maximum modulus principle and the Cauchy inequality.

\endsection

\section Index Theorem\par

\bigskip

In this section we prove that, under suitable hypotheses, the
section of \ref{siegellemma} has a small order of vanishing along
a point verifying the inequality of the main theorem. We will
prove an analogue of the ``Roth index theorem'' in this context.

\label vojta. remark\par\rmk In a first version of the paper we
deduced the index theorem from a generalization of the Vojta
version of Dyson lemma for curves [Vo]; but, due to the
``admissibility hypothesis'' in this kind of theorems, this could
be applied only in the case when {\it both} the $D_i$'s have the
same degree.\endrmk

Let ${\cal X}_1$ and ${\cal X}_2$ be the arithmetic surfaces. Let
$M$ be a line bundle over $({\cal X}_1 \times {\cal X}_2)_K$ and
$f \in H^0 (({\cal X}_1 \times {\cal X}_2)_K , M)$. We fix two
positive rational numbers $\vartheta_i \geq 1$.

Let $d_1$ and $d_2$ be two positive integers such that $d_i /
\vartheta_i \in {\bb Z}$.

Let $P := (P_1, P_2) \in ({\cal X}_1 \times {\cal X}_2)_K (K)$ be
a point and $z_i$ be local coordinate around $P_i$ in
$X_i:=(\X_i)_K$ ($i = 1, 2$). Let $e$ be a local generator of $M$
around $P$; consequently, near $P$, we can write $f = g \cdot e$
where $g$ is a regular function around $P$. We will say that $f$
has index at least $\delta$ in $P$ with respect to $d_1$ and $d_2$
and we will write ${\rm ind}_P (f, d_1, d_2) \geq \delta$ if, near
$P$, we write $g = \sum_{i,j} a_{i,j} z_1^i \cdot z_2^j$ and
$a_{i,j} = 0$ whenever
$$
{i \over d_1} \cdot \vartheta_1 + {j \over d_2} \cdot \vartheta_2
< \delta \, .
$$
The definition of the index is independent on the choices.

The condition ${\rm ind}_P(f, d_1, d_2)\geq\delta$ defines a
closed subscheme of $(\X_1\times\X_2)_K$ which will be denoted
$Z_\delta(f)$ (in the notation, the dependence on the $d_i$'s will
be clear from the context).

Let $M_i$ be the line bundles of generic degree one on ${\cal
X}_i$ ($i = 1, 2$) fixed in the previous section. As in the
previous section we will denote by $M_i$ the line bundle $pr_i^*
(M_i)$ on $X_1 \times X_2$ ($pr_i : {\cal X}_1 \times {\cal X}_2
\rightarrow {\cal X}_i$ being the natural projection).

The main theorem of this section is:

\label indextheorem. theorem\par\thm  Let $C$ and $\epsilon$ be
positive real numbers. Then we can find constants $B_j = B_j
(C,\epsilon)$ depending only on $C$, the $\vartheta_1$, and
$\epsilon$ (and on the other choices made until now), but
independent on the $d_i$'s, having the following property:

\noindent Suppose that:

\smallskip

\item{(a)} $f \in H^0 ({\cal X}_1 \times {\cal X}_2 ; M_1^{d_1} \otimes M_2^{d_2})$ is a
global section with ${\rm sup}_{\sigma \in M_{\infty}} \{ \Vert f \Vert_{\sigma} \} \leq C^{(d_1 + d_2)}$;

\smallskip

\item{(b)} the $d_i$'s are sufficiently big and divisible and $d_1 / d_2 \geq B_1$;

\smallskip

\item{(c)} $P := (P_1 , P_2) \in {\cal X}_1 \times {\cal X}_2 (K)$ is a rational point such that
$$
B_2 \leq h_{M_1} (P_1) \quad \hbox{and} \quad {h_{M_2} (P_2) \over
h_{M_1} (P_1)} \geq {d_1 \over d_2} \, ;
$$
then
$$
{\rm ind}_P (f , d_1 , d_2) \leq \epsilon \, .
$$
\endthm

\label faltingsproduct. remark\par\rmk The proof of the statement
above is directly inspired by the Faltings product theorem [Fa]
and can be deduced from it; we propose here a self contained proof
(which is simpler then the proof of the product theorem in this
situation).\endrmk

One can develop a height  for subvarieties of a fixed variety
(cf.[BGS]); this theory extends the height theory for points. We
will not recall here the definitions but we will recall the
properties of the heights that we need. Indeed, the only things we
need of the theories are the properties quoted below (consequently
a reader who do not know the theory can simply admit them).

We will use the following standard facts from the height theory of
subvarieties, one can find the proofs on [Fa2] or on [Ev]; if $Z$
is a closed subscheme of ${\cal X}_1 \times {\cal X}_2$ and $M$ is
an hermitian line bundle, then we denote by $h_M (Z)$ the height
of $Z$ with respect to $M$ as defined in [BGS]; by definition the
height of a closed subscheme is a real number. By linearity, the
height function is also defined on cycles:

\smallskip

\item{(a)} Suppose that $Z_i$ are closed irreducible reduced subschemes of ${\cal X}_i$ of relative dimension $\delta_i$ (over ${\bb Z}$) then
$$
h_{M_1^{d_1}\otimes M_2^{d_2}} (Z_1 \times Z_2) = (\delta_1 +
\delta_2 + 1)! \cdot d_1^{\delta_1} \cdot d_2^{\delta_2} \cdot
\left( {d_1 \cdot h_{M_1} (Z_1) \over (\delta_1 + 1)!} + {d_2
\cdot h_{M_2} (Z_2) \over (\delta_2 + 1)!} \right) \, ;
$$
this is proved in [Ev Lemma 8].

\smallskip

\item{(b)} Suppose that ${\cal X}_i = {\bb P}^1$ and $M_i = {\cal O} (1)$ and $C>1$ is a real
constant.  Then there is a constant $S$, depending only on $\X_i$
and the chosen metrics (but independent on the $d_i$'s and on
$C$), such that the following holds: let $f_1 , \ldots , f_r \in
H^0 ({\cal X}_1 \times {\cal X}_2 , {\cal O} (d_1) \otimes  {\cal
O} (d_2))$ be integral global sections such that ${\rm
sup}_{\sigma \in M_{\infty}} \{ \Vert f_i \Vert_{\sigma} \} \leq
C^{(d_1+d_2)}$; let $Y$ be the subscheme of ${\cal X}_1 \times
{\cal X}_2$ defined as the zero set of $\{ f_1 , \ldots , f_r \}$;
Let $X$ be an irreducible component of $Y$ with multiplicity $m_X$
then
$$
m_X \cdot h_{{\cal O} (d_1)\otimes{\cal O} (d_2)} (X) \leq
S\cdot\log(C) \cdot d_1 \cdot d_2 \cdot (d_1 + d_2) \, ;
$$
this is proved in [Fa2 Prop.2.17] or [Ev Lemma 9].

\smallskip

\item{(c)} If $f \in H^0({\cal X}_1 \times {\cal X}_2 , M_1^{d_1} \otimes  M_2^{d_2})$ then
$$
h_{M_1^{d_1} \otimes  M_2^{d_2}} ({\rm div} (f)) = h_{M_1^{d_1}
\otimes M_2^{d_2}} ({\cal X}_1 \times {\cal X}_2) + \sum_{\sigma
\in M_{\infty}} \int_{({\cal X}_1 \times {\cal X}_2)_{\sigma}}
\log \Vert f \Vert_{\sigma} (c_1 (M_1^{d_1} \otimes
M_2^{d_2})_{\sigma})^2 \, ;
$$
this is a direct consequence of the definition of height (cf.
[BGS]); consequently (using point (a)),  we can find a positive
constant $S$, depending only on the $\X_i$'s the $M_i$'s and the
chosen metrics, for which the following holds: let  $C>1$ be a
constant; if $f \in H^0({\cal X}_1 \times {\cal X}_2 , M_1^{d_1}
\otimes M_2^{d_2})$ is such that ${\rm sup}_{\sigma \in
M_{\infty}} \{ \Vert f \Vert_{\sigma} \} \leq C^{(d_1+d_2)}$ then
$$
h_{M_1^{d_1} \otimes  M_2^{d_2}} ({\rm div} (f)) \leq
S\cdot\log(C) \cdot d_1 \cdot d_2 \cdot (d_1 + d_2) \, .
$$

\medskip

\Proof {\it (of \ref{indextheorem}):} Let $f$ be the given section
and $Z$ be a geometrically irreducible reduced component of
$Z_{\epsilon} (f)$. Extending $K$ if necessary, we may suppose that
$Z$ is defined over $K$. It sufficies to prove that, under the
hypotheses of the theorem (with explicit and suitable $B_i$'s) the
point $P$ do not belong to $Z$. There are two cases, depending on
the dimension of $Z$.

\medskip

\noindent {\it Case 1: Dimension of $Z$ equal to one:} Let $Y :=
{\rm div} (f)$; since $Z$ is a divisor contained in $Y$ we have
$$
Y_K = m_Z \cdot Z + D \, ;
$$
where $D$ is an effective divisor on $({\cal X}_1 \times {\cal
X}_2)_K$. We claim that, if $d_1 / d_2 \geq \vartheta_1 /
(\epsilon \cdot \vartheta_2)$ then either there is a point $A \in
{\cal X}_2 (K)$ such that $Z = ({\cal X}_1)_K \times \{ A \}$, or
there is a point $B \in {\cal X}_1 (K)$ such that $Z = \{ B \}
\times ({\cal X}_2)_K$:

\label multiplicity. lemma\par\lemma Suppose that $Z$ is {\rm not}
as claimed,  then $m_Z \geq \epsilon \cdot {d_1 \over
\vartheta_1}$.
\endlemma

Let's show how the lemma implies the claim: Suppose that $Z$ is not
as claimed, then $(Z;M_1) > 0$; consequently, denoting by
$(\cdot;\cdot)$ the intersection pairing on the surface
$(\X_1\times\X_2)_K$,
$$
d_2 = (Y ; M_1) \geq \epsilon {{d_1}\over{\vartheta 1}} (Z ; M_1)
> \epsilon \cdot {d_1 \over \vartheta_1} \, ,
$$
so $d_1 / d_2 \leq \vartheta_1 / \epsilon$; thus, taking
$B_1>{{\vartheta_1}\over{\epsilon}}$, we find a  contradiction.

\bigskip

\Proof {\it (of the lemma):} Let $\eta$ be a generic point of $Z$
not contained in $D$; we may suppose that the restriction of both
projections are \'etale in a neighborhood of $\eta$. Let $z_1$ and
$z_2$ be local coordinates about the projections of $\eta$. In a
formal neighborhood of $\eta$, the divisor $Z$ is defined by a
irreducible element $h \in K [\![ z_1, z_2 ]\!]$ and $Y$ is
defined by the ideal $(h^{m_Z})$; because of our choice of $\eta$,
we have $h (z_1, z_2) = a_{10} z_1 + a_{01} z_2 + O((z_1 +
z_2)^2)$ with $a_{01} \cdot a_{10} \ne 0$, moreover, by definition
of $Z_{\epsilon} (f)$,
$$
(h (z_1 , z_2))^{m_Z} = \sum_{i,j} b_{ij} \cdot z_1^i \cdot z_2^j
$$
with $b_{ij} = 0$ whenever ${i \over d_1} \cdot \vartheta_1 + {j
\over d_2} \cdot \vartheta_2 \leq \epsilon$. Observe that
$b_{m_Z,0}=a_{01}^{m_Z}\neq 0$ thus
$$m_Z\geq \epsilon\cdot{{d_1}\over{\vartheta_1}}.$$

\smallskip

We thank the referee whose suggestions helped to drastically
simplify the proof of the lemma above.

\smallskip

We now come to the arithmetic part of the proof, in this case: $Z$
is either $({\cal X}_1)_K \times \{ A \}$ or $\{ B \} \times
({\cal X}_2)_K$ for suitable $A$ and $B$; remark that in the first
case $A=P_2$ and in the second case $B=P_1$. It is easy to see
that $m_Z$ is exactly  $\epsilon \cdot {d_2 \over \vartheta_2}$ in
the first case and exactly $\epsilon \cdot {d_1 \over
\vartheta_1}$ in the second: indeed it suffices to compute $m_Z$
on a smooth point of the support of $Z$ and $Y$. In the first
case, by applying properties (a) and (c) above, the fact that the
height is additive on cycles and the hypotheses, we can find an
explicit constant $R_1$ depending only on $C$ such that:
$$
\eqalign{ m_Z \cdot h_{M_1^{d_1} \otimes M_2^{d_2}} (Z) = & \ m_Z
\cdot d_1 (d_1 \cdot h_{M_1} ({\cal X}_1) + 2d_2 \cdot h_{M_2}
(A)) \cr \leq & \ h_{M_1^{d_1}\otimes M_2^{d_2}} ({\rm div} (f))
\leq R_1 \cdot d_1 \cdot d_2 (d_1 + d_2); \cr }
$$
consequently, since $m_Z = \epsilon \cdot {d_2 \over
\vartheta_2}$,
$$
{\epsilon \cdot d_1 \cdot d_2 \over \vartheta_2} \cdot (d_1 \cdot
h_{M_1} ({\cal X}_1) + 2d_2 \cdot h_{M_2} (A)) \leq R_1 \cdot d_1
\cdot d_2 (d_1 + d_2)
$$
thus
$$
d_2 \cdot h_{M_2} (A) < {R_2 \cdot \vartheta_2 \over \epsilon}
(d_1 + d_2) \, .
$$
Similarly, in the second case, we obtain
$$
{\epsilon \cdot d_1 \cdot d_2 \over \vartheta_1} \cdot (2 d_1\cdot
h_{M_1} (B) + d_2 \cdot h_{M_2} ({\cal X}_2)) \leq R_1 \cdot d_1
\cdot d_2 \cdot (d_1 + d_2) \, ,
$$
thus
$$
h_{M_1} (B) \leq {R_2 \over \epsilon} \cdot \left( 1 + {d_2 \over
d_1} \right) \, .
$$
This implies that, if $d_1 / d_2 \geq 1$, the point $(P_1, P_2)$
cannot be on $Z$ as soon as ${d_2 \over d_1} \cdot h_{M_2} (P_2)
\geq h_{M_1} (P_1) \geq {2R_2 \over \epsilon} \cdot {\rm max} \{
\vartheta_1 , \vartheta_2 \}$.

\bigskip

\noindent {\it Case 2: Dimension of Z equal to zero:} Denote by
$(P,Q) \in ({\cal X}_1 \times {\cal X}_2)_K (K)$ the support of
$Z$. In this case we need to project on ${\bb P}^1 \times {\bb
P}^1$. We fix once for all a finite set of coverings $\gamma_{i,j}
: ({\cal X}_i)_K \rightarrow {\bb P}^1$ with the following
property: if $U_{i,j} \subseteq ({\cal X}_i)_K$ is the open set
over which $\gamma_{i,j}$ is \'etale, then $\bigcup_j U_{i,j} =
({\cal X}_i)_K$ and $\gamma_{i,j}^* ({\cal O} (1)) \simeq
(M_i)_K^{t_i}$ for suitable $t_i$ (we fix such isomorphisms). We
also suppose that each $\gamma_{i,j}$ extends to a generically
finite morphism $\gamma_{i,j} : {\cal X}_i \rightarrow {\bb
P}_{O_K}^1$ (this can be obtained after a suitable blow up of
${\cal X}_i$). We equip the line bundle ${\cal O} (1)$ on ${\bb
P}^1$ with the Fubini-Study metric $\Vert \cdot \Vert_{FS}$. Fix a
constant $A$ such that
$$
A^{-1} \gamma_{i,j}^\ast (\Vert \cdot \Vert_{FS}) \leq \Vert \cdot
\Vert_{M_i}^{t_i} \leq A \gamma_{i,j}^\ast (\Vert \cdot
\Vert_{FS}) \, .
$$

We may suppose that $(P,Q) \in ({\cal X}_1 \times {\cal X}_2)(K)$
is contained in $U_{1,1} \times U_{2,1}$. Denote by $\Gamma$ the
morphism $\gamma_{1,1} \times \gamma_{1,2} : {\cal X}_1 \times
{\cal X}_2 \rightarrow {\bb P}^1 \times {\bb P}^1$. Put $d_i = t_i
\cdot a_i$; then $\Gamma^* ({\cal O} (a_1 , a_2)) =
M_1^{d_1}\otimes M_2^{d_2}$ and $g := \Gamma_* (f) \in H^0 ({\bb
P}^1 \times {\bb P}^1 , {\cal O} (d_1 \cdot t_2 , d_2 \cdot t_1)$.
It is easy to verify that there exists an absolute constant $A_1$
such that
$$
\Vert g \Vert_{FS} \leq A_1^{(d_1 + d_2)} \Vert f \Vert
$$
and that $(P' , Q') := \Gamma (P,Q)$ is contained in $Z_{\epsilon}
(g)$. Consequently it suffices to prove the theorem when ${\cal
X}_1 = {\cal X}_2 = {\bb P}^1$, $M_1 = M_2 = {\cal O} (1)$ and
$\O(1)$ is equipped with the Fubini-Study metric.

We first look to the irreducible components $Z'$ of $Z_{\epsilon /
2}$ containing $(P' , Q')$. If there is such a $Z'$ of dimension
one, then we are reduced to the previous case and we are done. We
may then suppose that the support of $Z'$ is $(P' , Q')$ too. Let
$I_{\epsilon}$ and $I_{\epsilon / 2}$ be the ideal of $Z$ and $Z'$
in the completion $K [\![ z_1 , z_2 ]\!]$ of the local ring of
${\bb P}^1 \times {\bb P}^1$ in $(P' , Q')$; let $h = \alpha \cdot
z_1^{r_1} z_2^{r_2} + \ldots$ be an element of $I_{\epsilon / 2}$
then
$$
{\partial^{i_1 + i_2} \over \partial z_1^{i_1} \cdot \partial
z_2^{i_2}} h = \alpha_1 z_1^{(r_1 - i_1)'} z_2^{(r_2 - i_2)'} +
\ldots
$$
(where $(a)' := {\rm sup} \{ a , 0\}$) for a suitable $\alpha_1$;
and $\alpha_1$ is zero only if $\alpha$ is zero or $\alpha \ne 0$
and one of the $(r_j - i_j)'$ is zero. If ${i_1 \over d_1} \cdot
\vartheta_1 + {i_2 \over d_2} \cdot \vartheta_2 < {\epsilon \over
2}$ and $h \in I_{\epsilon / 2}$ then ${\partial^{i_1 + i_2} \over
\partial z_1^{i_1} \cdot \partial z_2^{i_2}} h \in I_{\epsilon}
\subseteq (z_1 , z_2)$. This implies that $h$, and consequently
$I_{\epsilon / 2}$, is contained in the ideal $(z_1^{\epsilon d_1
/ (4 \cdot \vartheta_1)} , z_2^{\epsilon d_2 / (4\vartheta_2)})$.
Thus, the multiplicity of $Z'$ in $Z_{\epsilon / 2} (g)$ is at
least ${1 \over \vartheta_1 \cdot \vartheta_2} \cdot \left(
{\epsilon \over 4} \right)^2 \cdot d_1 d_2$.

Every differential operator ${\partial^{i_1 + i_2} \over \partial
z_1^{i_1} \partial z_2^{i_2}}$ with ${i_1 \over d_1} \cdot
\vartheta_1 + {i_2 \over d_2} \cdot \vartheta_2 \leq {\epsilon
\over 2}$ can be seen as a linear endomorphism $D^{(i_1,i_2)}$ of
$H^0({\bb P}^1 \times {\bb P}^1 , {\cal O}(d_1, d_2))$. For every
infinite place $\sigma \in M_{\infty}$ the norm of the operator
$D^{(i_1,i_2)}$ (with $i_1$ and $i_2$ bounded as above) is bounded
from above by $2^{{\rm max} \{\vartheta_1 , \vartheta_2 \} \cdot
(d_1+ d_2)}$. We apply property (b) above and the hypotheses and
we find a constant $R'$, depending only on $C$
$$
m_{Z'} \cdot (d_1 \cdot h_{{\cal O} (1)} (P') + d_2 \cdot h_{{\cal
O} (1)} (Q')) \leq R' \cdot d_1 \cdot d_2 \cdot (d_1 + d_2);
$$
consequently, since $m_{Z'} \geq {1 \over \vartheta_1 \cdot
\vartheta_2} \cdot \left( {\epsilon \over 4} \right)^2 \cdot d_1
\cdot d_2$, we obtain
$$
d_1 \cdot h_{{\cal O} (1)} (P') + d_2 \cdot h_{{\cal O} (1)} (Q')
\leq \vartheta_1 \cdot \vartheta_2 \cdot \left( {4 \over \epsilon}
\right)^2 \cdot R' (d_1 + d_2) \, .
$$
If we suppose that $\vartheta_1 \cdot \vartheta_2 \cdot \left( {4
\over \epsilon} \right)^2 \cdot R' \leq h_{{\cal O} (1)} (P_1)
\leq h_{{\cal O} (1)} (P_2)$ the point $P$ cannot belong to $Z'$
and this concludes the proof of the theorem.

\label variationofBis. rmk\par\rmk We observe that, from the proof
one deduce that the constants $B_i$'s of theorem
\ref{indextheorem} may be chosen of the form $B_i=S_i\log(C)$,
where the $S_i$ depend only on the $\X_j$ the $M_i$ and the chosen
metrics (but independent on $C$).

\endrmk

\endsection

\section Generalized Cauchy inequalities\par

\bigskip

Fix $\vartheta_i \in {\bb Q}_{\geq 1}$ and the divisors $D_i : {\rm
Spec} (O_L) \rightarrow {\cal X}_i$ as in section 3. For every
rational positive $\delta$ and couple of positive integers $(d_1,
d_2)$, let ${\cal I}_{\underline{\vartheta} , \delta ,
\underline{d}}$ be the ideal sheaf of ${\cal X}_1 \times {\cal X}_2$
defined in section 3. Let $p : \tilde {\cal X}_{\delta} \rightarrow
{\cal X}_1 \times {\cal X}_2$ be the blow up along ${\cal
I}_{\underline{\vartheta} , \delta , \underline{d}}$ and let
$E_{\delta}$ be the corresponding exceptional divisor on it. We can
find a very small positive constant $\alpha$ such that, if the $d_i$
are sufficiently big, there is a surjection
$$
\beta_{\delta}\colon\bigoplus_{\delta \leq {i_1 \over d_1} \cdot
\vartheta_1 + {i_2 \over d_2} \cdot \vartheta_2 \leq \delta +
\alpha} {\cal O} (-i_1 \cdot D_1) \otimes {\cal O} (-i_2 \cdot D_2)
\longrightarrow\!\!\!\!\rightarrow {\cal I}_{\underline{\vartheta} ,
\delta , \underline{d}} \, .
$$
Observe that $\alpha$ is independent on the $d_i$'s, provided that
they are sufficiently big.

To simplify notations we will denote by $H$ the set
$$
\left\{ (i_1 , i_2) \in {\bb Z} \times {\bb Z}\;\; \vert\;\;
\delta \leq {i_1 \over d_1} \cdot \vartheta_1 + {i_2 \over d_2}
\cdot \vartheta_2 \leq \delta + \alpha \right\} \, .
$$
If $M$ is an hermitian line bundle on ${\cal X}_1 \times {\cal
X}_2$, by abuse of notation, we will denote again by $M$ the pull
back of $M$ to $\tilde{\cal X}_{\delta}$.

The surjection $\beta_{\delta}$ above induces a surjection
$$
\beta_{\delta}\colon\bigoplus_{(i_1,i_2)\in H}  {\cal O} (-i_1 \cdot
D_1) \otimes {\cal O} (-i_2 \cdot D_2)
\longrightarrow\!\!\!\!\rightarrow \O_{\tilde
X_{\sigma}}(-E_{\delta});$$ consequently the line bundle ${\cal
O}_{\tilde X_{\sigma}} (E_{\delta})$ is naturally equipped with the
structure of {\it hermitian line bundle}.

If $P_i\in\X_i(K)$ are $K$-rational points of $\X_i$, they extend to
sections $P_i\colon B:=\Spec(O_K)\to\X_i$. We will denote by
$P\colon B\to \X_1\times\X_2$ the section $P_1\times P_2$ and by
$\tilde P\colon B\to \tilde\X_{\delta}$ the strict transform of $P$.

The theorem we want to prove in this section is the following:

\label gencauchyineq. theorem\par\thm Let $M$ be an hermitian line
bundle on $\X_1\times\X_2$ and $A$ and $\epsilon$ be  positive
constants. There is a constant $C$ depending only $A$, on the
models, the metrics, the $\vartheta_i$'s etc. but independent on the
$d_i$'s for which the following holds:

\noindent Let  $f \in H^0 ({\cal X}_1 \times {\cal X}_2 , M \otimes
{\cal I}_{\underline{\vartheta} , \delta , \underline{d}})$ be a
global section such that ${\rm sup}_{\sigma \in M_{\infty}} \{ \Vert
f \Vert_{\sigma} \} \leq A$. Let  $P := P_1 \times P_2 : {\rm Spec}
(O_K) \rightarrow {\cal X}_1 \times {\cal X}_2$ be a rational point
such that ${\rm ind}_P (f, d_1, d_2) \leq \epsilon$; then there
exists $\epsilon' \leq \epsilon$, two positive integers $i_1$ and
$i_2$ such that ${i_1 \over d_1} \cdot \vartheta_1 + {i_2 \over d_2}
\cdot \vartheta_2 \leq \epsilon$ and a {\rm non zero} global section
$\tilde f \in H^0 (\tilde P , M \otimes \omega_{{\cal X}_1 /
B}^{i_1} \otimes \omega_{{\cal X}_2 / B}^{i_2} \otimes {\cal O}
(-E_{\delta - \epsilon'}))$ such that
$$
\sup_{\sigma \in M_{\infty}} \{ \Vert \tilde f \Vert_{\sigma} \}
\leq A \cdot C^{(d_1 + d_2)} \, .
$$
\endthm

Before we start the proof of the theorem, we need to introduce
some notations and some tools.

Let ${\cal L}_1$ and ${\cal L}_2$ be two line bundles on ${\rm
Spec} (O_K)$. For every couple of positive integers $(i_1, i_2)$
we define the differential operator
$$
D^{(i_1 , i_2)} : {\cal O} [\![ {\cal L}_1 , {\cal L}_2 ]\!]
\longrightarrow {\cal O} [\![ {\cal L}_1 , {\cal L}_2 ]\!] \otimes
{\cal L}_1^{i_1} \otimes {\cal L}_2^{i_2}
$$
in the following way: let $e_1$ (resp. $e_2$) be a local generator
of ${\cal L}_1$ (resp. of ${\cal L}_2$) then we define
$$
D^{(i_1 , i_2)} (e_1^a \otimes e_2^b) :=\left\{\eqalign{ \left( {a
\atop i_1} \right) \cdot \left( {b \atop i_2} \right) \cdot
e_1^{a-i_1} \otimes e_2^{b - i_2} \otimes (e_1^{i_1} \otimes
e_2^{i_2}) & \;\; {\rm if} \; a\geq i_1 \; {\rm and}\; b\geq
i_2\cr 0\; \;\; \;\; \; \; \;\; \;\; \;\; \;\; \;\; \;&\; \; {\rm
otherwise;}\cr}\right.
$$
and extend it linearly to ${\cal O} [\![ {\cal L}_1, {\cal L}_2
]\!]$; one easily verify that this definition do not depends on the
choice of the local generators. The module ${\cal O} [\![ {\cal
L}_1, {\cal L}_2 ]\!] \otimes {\cal L}_1^{i_1} \otimes {\cal
L}_2^{i_2}$ has a natural structure of ${\cal O} [\![ {\cal L}_1,
{\cal L}_2 ]\!]$-module (multiplication on the right). one can
easily verify that $D^{(i_1, i_2)}$ is a differential operator: it
is $O_K$-linear (by definition) and it satisfy the (iterated)
Leibnitz-rule; for instance $D^{n,0} (f \cdot g) = \sum \left( {n
\atop i} \right) \cdot D^{(i,0)} (f) \cdot D^{(n-i,0)} (g)$
$(D^{(i,0)} (f) \in {\cal O} [\![ {\cal L}_1 , {\cal L}_2 ]\!]
\otimes {\cal L}_1^i$ and $D^{(n-i,0)} (g) \in {\cal O} [\![ {\cal
L}_1 , {\cal L}_2 ]\!] \otimes {\cal L}_1^{n-i})$, consequently
$D^{(i,0)} (f) \cdot D^{(n-i,0)} (g) \in {\cal O} [\![ {\cal L}_1 ,
{\cal L}_2 ]\!] \otimes {\cal L}_1^n)$.

If $\sigma \in M_{\infty}$ is an infinite place, then ${\cal O}
[\![ {\cal L}_1, {\cal L}_2 ]\!]_{\sigma}$ is (non canonically)
isomorphic to the ring of formal power series in two variables and
the operators $D^{(a,b)}$ are the usual iterated derivatives.

Although it is not necessary, we will tacitly authorize ourself to
pass to the Hilbert class field extension: consequently we will
suppose that every line bundle on $B$ is trivial; this is not
necessary, but highly simplify the notations.

Denote by $(\widehat{{\cal X}_1 \times {\cal X}_2})_P$ the formal
completion of ${\cal X}_1 \times {\cal X}_2$ around $P$. By
\ref{completionontwopoints}, we find a canonical isomorphism
$$
\Psi_P : (\widehat{{\cal X}_1 \times {\cal X}_2})_P
\buildrel{\sim}\over{\longrightarrow} {\rm Spf} ({\cal O} [\![
{\cal O} (-P_1) \vert_{P_1} , {\cal O} (-P_2) \vert_{P_2} ]\!] )
\, .
$$
We will denote by $I_P \subset {\cal O} [\![ {\cal O}
(-P_1)\vert_{P_1} , {\cal O} (-P_2)\vert_{P_2} ]\!]$ the ideal
corresponding to the ideal of definition of $(\widehat{{\cal X}_1
\times {\cal X}_2})_P$ defining the point section $P$ (with the
reduced structure).

\bigskip

\Proof ({\it Of theorem 5.1}) Let $p_i : {\cal X}_1 \times {\cal
X}_2 \rightarrow {\cal X}_i$ the projection. Denote by $\I_{D_i}$
the restriction of the ideal sheaf $p_i^* ({\cal O} (-D_i))$ to
$(\widehat{{\cal X}_1 \times {\cal X}_2})_P$. The image of
$\I_{D_i}$ by $\Psi$ is a principal ideal of ${\cal O} [\![ {\cal
O}(-P_1)\vert_{P_1} , {\cal O}(-P_2)\vert_{P_2} ]\!]$ generated by
an element $G_i$. If $\delta$ is a positive rational number, we
denote then by $I_{\delta , \underline{d}} \subset {\cal O} [\![
{\cal O} (-P_1) \vert_{P_1} , {\cal O} (-P_2) \vert_{P_2} ]\!]$
the ideal generated by the elements $G_1^i \cdot G_2^i$ with ${i
\over d_1} \cdot \vartheta_1 + {j \over d_2} \cdot \vartheta_2
\geq \delta$. The ideal  $I_{\delta , \underline{d}}$ is the
image, via $\Psi$, of the restriction to $(\widehat{{\cal X}_1
\times {\cal X}_2})_P$ of the ideal sheaf ${\cal
I}_{\underline{\vartheta} , \delta , \underline{d}}$.
Consequently, a global section $f \in H^0 ({\cal X}_1 \times {\cal
X}_2 ; M \otimes {\cal I}_{\underline{\vartheta} , \delta ,
\underline{d}})$ restricted to $(\widehat{{\cal X}_1 \times {\cal
X}_2})_P$ will determine an element
$$
F = \sum_{{i \over d_1} \cdot \vartheta_1 + {j \over d_2} \cdot
\vartheta_2 \geq \delta} a_{ij} \cdot G_1^i \cdot G_2^j \, .
$$
If $(i_1, i_2)$ is a couple of indices such that ${i_1 \over d_1}
\cdot \vartheta_1 + {i_2 \over d_2} \cdot \vartheta_2 \leq
\epsilon$ then a direct computation using the iterated Leibnitz
rule gives $D^{(i_1,i_2)} (F) \in I_{\delta - \epsilon ,
\underline{d}} \otimes M \vert_P \otimes {\cal O} (-i_1 P_1)
\vert_{P_1} \otimes {\cal O} (-i_2 P_2) \vert_{P_2}$.

Since the  index ${\rm ind}_P (f, d_1, d_2)$ of $f$ at $P$ is less
or equal then $\epsilon$, then we can find a couple of positive
integers $(i_1, i_2)$ such that ${i_1 \over d_1} \cdot \vartheta_1 +
{i_2 \over d_2} \cdot \vartheta_2 \leq \epsilon$ and such that the
class $\tilde f$ of $D^{(i_1 , i_2)} (f)$ in $({\cal O} [\![ {\cal
O} (-P_1) \vert_{P_1} , {\cal O} (-P_2) \vert_{P_2} ]\!] \otimes M
\vert_P \otimes {\cal O} (-i_1 P_1) \vert_{P_1} \otimes {\cal O}
(-i_2 P_2) \vert_{P_2}) / I_P \simeq H^0 (P,M \otimes {\cal O} (-i_1
P_1) \vert_{P_1} \otimes {\cal O} (-i_2 P_2) \vert_{P_2})$ is {\it
non zero}. Thus, using adjunction formula, we find a {\it non zero}
section in $\tilde f \in H^0 (P, M \otimes \omega_{{\cal X}_1 /
B}^{i_1} \vert_{P_1} \otimes \omega_{{\cal X}_2 / B}^{i_2}
\vert_{P_2} \otimes {\cal I}_{\underline{\vartheta} ,
\delta-\epsilon , \underline d})$.

Let $\sigma \in M_{\infty}$ be an infinite place. we fix once for
all a covering of $({\cal X}_i)_{\sigma}$ by open sets $U_{ij}$
analytically equivalent to a disk (with coordinate $z$)for which the
following holds:

\noindent -- The line bundle $\O(\Delta_i)$ is trivial on
$U_{ij}\times U_{ik}$; and we fix once for all a trivialization.

\noindent -- The line bundle $M_{\sigma}$ is trivial on
$U_{1,j}\times U_{2,k}$.

Let $\Vert\cdot\Vert_\ell$ be the metric on the line bundle
$\O(P_\ell)_{\sigma}$. Let $\Bbb I_\ell$ be the canonical section of
$\O(P_\ell)_{\sigma}$. There is a ${\cal C}^{\infty}$ function
$\rho_{\ell ij}$ on $U_{ij}$ such that

$$\Vert \Bbb I_\ell\Vert_{\ell}(z)=\rho_{\ell, ij}(z)\cdot\vert z-z(P_\ell)\vert.$$

Due to our choices, we can find (and fix once for all) two constants
$A_1$ and $A_2$ {\it independent on the $P_i$} such that
$$A_1\leq \rho_{\ell, ij}(z)\leq A_2.$$
Thus, we apply \ref{cauchylemma} and we find an absolute constant
$C_1$, independent on $P$ and on the $d_i$, such that

$$
\sup \{ \Vert \tilde f \Vert_{\sigma} \} \leq A \cdot C_1^{(d_1 +
d_2)} \, .
$$

The section $\tilde f$ extends to a section, denote it again by
$\tilde f$, of $(M \otimes \, \omega_{{\cal X}_1}^{i_1} \otimes \,
\omega_{{\cal X}_2}^{i_2})_{\sigma}$ on a neighborhood of $P$, which
we may suppose to be one of the  products of the $U_i$ above; a
similar argument shows that $\sup \{ \Vert \tilde f \Vert \} \leq A
\cdot C_1^{(d_1 + d_2)}$.

Let $\tilde{\cal X} \rightarrow {\cal X}_1 \times {\cal X}_2$ be the
blow up along the ideal ${\cal I}_{\underline{\vartheta} , \delta -
\epsilon , \underline{d}}$ and $E_{\delta - \epsilon}$ be the
exceptional divisor; let $\tilde P : {\rm Spec} (O_K) \rightarrow
\tilde X$ be the strict transform of $P$. By definition $\tilde f$
will give a non zero section (which we will denote with the same
symbol) $\tilde f \in H^0 (\tilde P , M \otimes \omega_{{\cal X}_1 /
B}^{i_1} \otimes \omega_{{\cal X}_2 / B}^{i_2} (-E_{\delta -
\epsilon}))$. We will now give an upper bound for the norm of
$\tilde f$. As before, once we take a suitably chosen (once for all)
open covering of $({\cal X}_i)_{\sigma}$, in the analytic topology,
the existence of the upper bound as in the statement of the theorem
is consequence of \ref{localcauchy} below.
\smallskip
Let ${\bb D}$ be an open disk, $0 \in {\bb D}$ be a point on it
and $z$ be a coordinate with a simple zero on $0$. Suppose that
$\rho_i (z)$ ($i = 1, 2$) are two ${\cal C}^{\infty}$ functions on
${\bb D}$; suppose that we can find two positive constants $B_1$
and $B_2$ such that $B_1 \leq \rho_i (z) \leq B_2$. We define two
metrics $\Vert \cdot \Vert_i$ on ${\cal O} (0)$ by the formula
$\Vert {\bb I}_0 \Vert_i = \vert z \vert \rho_i (z)$.

Let $p_i : {\bb D} \times {\bb D} \rightarrow {\bb D}$ the $i$-th
projection, we will denote by ${\cal O} (-0_i)$ the line bundle
$p_i^* ({\cal O} (-0))$ and by $z_i$ the holomorphic function $p_i
(z)$ (it is the canonical section of ${\cal O} (0_i)$). We will
suppose that ${\cal O} (0_i)$ is equipped with the pull-back, via
$p_i$ of the metric $\Vert \cdot \Vert_i$.

Fix positive rational numbers $\vartheta_i$ and $\delta$. For
every couple of sufficiently divisible positive integers $(d_1,
d_2)$ define ${\cal I}_{\underline{\vartheta} , \delta ,
\underline{d}}$ to be the ideal sheaf of ${\cal O}_{{\bb D} \times
{\bb D}}$ generated by the monomials $z_1^{i_1} \cdot z_2^{i_2}$
with ${i_1 \over d_1} \cdot \vartheta_1 + {i_2 \over d_2} \cdot
\vartheta_2 \geq \delta$.

Let $b : \tilde X \rightarrow {\bb D} \times {\bb D}$ be the blow up
of ${\cal I}_{\underline{\vartheta} , \delta , \underline{d}}$ and
$E:=E_{\delta} \subset \tilde X$ be the exceptional divisor. In the
same way as before, if the $d_i$ are sufficiently big, we have a
surjection
$$
\bigoplus_{(i_1,i_2)\in H} {\cal O} (-i_1 \cdot 0_1) \otimes {\cal
O} (-i_2 \cdot 0_2) \longrightarrow\!\!\!\!\rightarrow {\cal
I}_{\underline{\vartheta} , \delta , \underline{d}} \, .
$$
which induces a metric on $\O(E)$.

\label localcauchy. theorem\par\thm There exists a constant $B$
depending only on $\vartheta_i , \delta$ and the constants $A_i$
such that if the $d_i$'s are sufficiently big and divisible, $f
\in H^0 ({\bb D} \times {\bb D} , {\cal I}_{\underline{\vartheta}
, \delta , \underline{d}})$ and $\tilde f$ is the corresponding
section in $H^0(\tilde X , {\cal O}(-E))$ then, for every
$z\in\tilde X$,
$$
\Vert \tilde f \Vert (z)\leq \Vert f \Vert(b(z)) \cdot B^{(d_1 +
d_2)} \, .
$$
\endthm

\Proof Denoting by ${\bf P}$ the projective bundle ${\bf Proj}
(\bigoplus_{(i_1 , i_2) \in H} \, {\cal O} (-i_1 \cdot 0_1) \otimes
{\cal O} (-i_2 \cdot 0_2))$ over $\bb D \times \bb D$ we get a
commutative diagram
$$
\matrix{ \tilde X &\buildrel{\iota}\over{\hookrightarrow} &{\bf P}
\cr &\searrow &\downarrow \cr && \ \ {\bb D} \times {\bb D} \, .
\cr }
$$
Moreover, by construction we have an {\it isometry} $\iota^* ({\cal
O} (1)) \simeq {\cal O} (-E)$. Remark that ${\bf P}$ is isomorphic
to ${\bb D} \times {\bb D} \times {\bb P}^N$ for a suitable $N$.
Denote by $[u_{i_1 , i_2}]_{(i_1,i_2)\in H}$  the homogeneous
coordinates on ${\bb P}^N$; the blow up $\tilde X$ is defined by the
equations
$$
u_{j_1 , j_2} \cdot z_1^{i_1} \cdot z_2^{i_2} = u_{i_1 , i_2}
\cdot z_1^{j_1} \cdot z_2^{j_2}
$$
for all $(i_1, i_2)$ and $(j_1, j_2)$ in $H$.

Let's work on the local chart $u_{i_1 , i_2} \ne 0$; a local
computation shows that over this chart \labelf normofE\par$$ \Vert
E \Vert = {\vert z_1^{i_1} \cdot z_2^{i_2} \vert \over \vert
u_{i_1 , i_2} \vert} \cdot \sqrt{\sum_{(j_1 , j_2) \in H} \left(
\vert u_{j_1 , j_2} \vert \cdot \rho_1^{j_1} \cdot \rho_2^{j_2}
\right)^2} \, .\eqno {{(\numfo)}}
$$\advance\ssnu by1
Let $f \in H^0 ({\bb D} \times {\bb D}, I_{\underline{\vartheta} ,
\delta , \underline{d}})$. The pull-back $b^* (f)$ naturally
defines a global section $\tilde f \in H^0 (\tilde X , {\cal O}
(-E))$. Over the chart $u_{i_1 , i_2} \ne 0$ we can find a
holomorphic function $h$ such that $f = z_1^{i_1} \cdot z_2^{i_2}
\cdot h$. In order to conclude the proof of the theorem we have to
give an upper bound for \labelf normofh\par$$ \vert h \vert \cdot
{\sqrt{\sum \vert u_{j_1 , j_2} \vert^2 \cdot \rho_1^{2j_1} \cdot
\rho_2^{2j_2}} \over \vert u_{i_1 , i_2} \vert} \, . \eqno
{{(\numfo)}}
$$
Fix a very small positive $\epsilon$; we may suppose that we
are in the disk
$$
{\vert u_{j_{1} , j_{2}} \vert \over \vert u_{i_1 , i_2} \vert}
\leq 1 + \epsilon \, ;
$$
if this is not verified, it suffices to change the local chart.
consequently, we can find a constant $B_1$ depending only on the
norms (in particular independent on the $d_i$'s) for which  the
expression in \ref{normofh} is bounded from above by
$$
\vert h \vert \cdot B_1^{(d_1 + d_2)} \, .
$$
Since $h$ is holomorphic, the function $\vert h \vert$ will assume
its maximum on the border. We may assume that the $d_i$ are such
that ${d_i \cdot \delta \over \vartheta_i} \in {\bb N}$. On our
chart, $z_1^{{\delta \cdot d_1 \over \vartheta_1}} = z_1^{i_1} \cdot
z_2^{i_2} \cdot u_{{\delta \cdot d_1 \over \vartheta_1} , 0}$ (resp.
$z_2^{{\delta \cdot d_2 \over \vartheta_2}} = z_1^{i_1} \cdot
z_2^{i_2} \cdot u_{0 , {\delta \cdot d_1 \over \vartheta_1}}$) and
$\vert u_{{\delta \cdot d_2 \over \vartheta_2} , 0} \vert$ (resp.
$u_{0 , {\delta \cdot d_2 \over \vartheta_2}}$) is less or equal to
$1 + \epsilon$. Consequently, if $\vert z_1^{{\delta \cdot d_1 \over
\vartheta_1}} \vert = 1$ (resp. $\vert z_2^{{\delta \cdot d_2 \over
\vartheta_2}}\vert = 1$) then $1 \leq \vert z_1^{i_1} \cdot
z_2^{i_2} \vert \cdot (1 + \epsilon)$ thus
$$
\vert h \vert \leq {\Vert f \Vert \over \vert z_1^{i_1} \cdot
z_2^{i_2} \vert} \leq (1+\varepsilon) \cdot \Vert f \Vert \, ;
$$
the conclusion of the theorem easily follows.

\endsection

\section Proof of the main theorem\par

\bigskip

In this section we will give the proof of the main theorem of the
paper: Theorem \ref{maindyson}.

We recall all the tools and the ingredients: ${\cal X}_i$ are two
regular arithmetic surfaces projective over $B := {\rm Spec}
(O_K)$ over which we fixed arithmetically ample hermitian line
bundles $M_i$ and symmetric hermitian metrics on ${\cal O}
(\Delta_i)$ ($\Delta_i$ being the diagonal on ${\cal X}_i \times
{\cal X}_i$). Eventually we fix a place $\sigma \in M_K$.

We fix two finite extensions $L_i$ of $K$ and two reduced divisors
$D_i : B_{L_i} := {\rm Spec} (O_{L_i} ) \rightarrow {\cal X}_i$.
We denote by $L$ the composite field $L_1 \cdot L_2$ and by $n$
the degree of the extension $L / K$. We fix two positive rational
numbers $\vartheta_i \geq 1$ and a positive $\epsilon$ such that
$\vartheta_1 \cdot \vartheta_2 \geq 2n + \epsilon$. We will denote
by $T(D_i)$ the positive real number introduced in \S 2. We fix a
function $\varphi: S\to [0;1]$ such that $\sum_{v\in
S}\varphi(v)=1$. Theorem 2.1 will be consequence of the following:

\bigskip

\label proofofmaindyson. theorem\par\thm There exists a constant
$A$ depending only on the arithmetic surfaces ${\cal X}_i$, the
$\vartheta_i$, the $\epsilon$, the hermitian line bundles $M_i$,
the symmetric metrics on the diagonals ${\cal O} (\Delta_i)$, the
set $S$ and the function $\varphi$, for which the following holds:

Let $D_i \subset {\cal X}_i$ be divisors as above, and $P_i \in
{\cal X}_i (B)$ be two rational sections such that

\smallskip

\item{(i)} For every place $v\in S$, we have that $\lambda_{D_1 , v} (P_1) >\varphi(v)\cdot\vartheta_1 \cdot h_{M_1} (P_1)$
and $\lambda_{D_2 , v} (P_2) > \varphi(v)\cdot\vartheta_1 \cdot
h_{M_2} (P_2)$;

\smallskip

\item{(ii)} $h_{M_1} (P_1) \geq A \cdot T (D_1) \cdot T (D_2)$.

\smallskip

Then
$$
h_{M_2} (P_2) \leq A \cdot T (D_1) \cdot T(D_2) \cdot h_{M_1}
(P_1) \, .
$$
\endthm
\medskip

\Proof We first treat the case when, for at least one place, each
of the $P_i$'s is "far from $D_i$". Suppose that $v\in S$ is an
infinite place, then take a covering of $({\cal X}_i)_{v}$ by open
sets $U_{ij}$, analytically equivalent to the disk of radius $1$
and such that the open subsets analytically equivalent to the disk
of radius $1/2$ also cover the $({\cal X}_i)_{v}$. We can then
find a constant $A_2$ such that if $U_{ij_k}$ are the open sets
containing the $(D_i)_{v}$ and $(P_i)_{v}$ are not contained in
the $U_{ij_k}$ then $\lambda_{D_i , v} (P_i) \leq A_2$.
Consequently, we see that, taking $A$ much bigger then $A_2$
(which is independent on the $D_i$), in this case condition (i)
and condition (ii) are in contradiction. In particular the theorem
holds in this case. A similar argument holds if $v$ is a finite
place.

Suppose that $\vartheta_1 \cdot \vartheta_2 = 2n + \epsilon$;
define $\epsilon_1 := {\epsilon \over n+1}$ and $\delta := 2 +
\epsilon_0$. a suitable choice of $\epsilon_0$ allows to suppose
that the hypotheses of theorem \ref{smallsection} are verified.

Here again, "absolute constant" will be equivalent to say "a
constant which depends only on the $\X_i$, the hermitian line
bundles $M_i$, the metrics on the diagonals, the $\vartheta_i$'s,
but independent on the $D_i$'s and on the $d_i$'s".

For every couple of positive integers $d_1$ and $d_2$, let ${\cal
I}_{\underline{\vartheta} , \delta , \underline d}$ be the ideal
sheaf on ${\cal X}_1 \times {\cal X}_2$ defined in \S3 and having
support on $D_1 \times D_2 \subset {\cal X}_1 \times {\cal X}_2$.

Fix an absolute constant $A_3$ such that $h_{\omega_{{\cal X}_i /
B}} (\cdot) \leq A_3 \cdot h_{M_i} (\cdot)$ and let $\epsilon_2$
such that $\epsilon_2 < {\epsilon_1 \over 1+2 \cdot A_3}$.

We apply \ref{smallsection} and we find an absolute constant
constant $A_4$ such that, each time $d_i$'s are sufficiently big and
divisible we can find a non zero global section $f \in H^0 ({\cal
X}_1 \times {\cal X}_2 , M_1^{d_1} \otimes M_2^{d_2} \otimes {\cal
I}_{\underline{\vartheta} , \delta , \underline d})$ such that
$$
\sup_{\sigma \in M_{\infty}} \{ \log \Vert f \Vert_{\sigma} \} \leq
A_4 \cdot T(D_1) \cdot T(D_2) (d_1 + d_2) \, .
$$

One apply Theorem \ref{indextheorem} with $\log(C) = A_3 \cdot
T(D_1) \cdot T(D_2)$ and $\epsilon = \epsilon_2$ and deduce the
existence of a constant $A_5$ for which, if a point $P$ verify
(a), (b) and (c) of loc cit. then the index
$ind_P(f,d_1,d_2)<\epsilon_2$; by remark \ref{variationofBis} one
can see that $A_5$ is again of the form $A_6 \cdot T(D_1) \cdot
T(D_2)$ with $A_6$ independent on the $D_i$'s.

Suppose that $P_i : B \rightarrow {\cal X}_i$ are two sections
which satisfy hypothesis (i) and such that
$$
h_{M_2} (P_2) > A_6 \cdot T(D_1) \cdot T(D_2) \, h_{M_1} (P_1)
$$
we will prove that there exists a constant $A_7$ such that $h_{M_1}
(P_1) \leq A_7 \cdot T(D_1) \cdot T(D_2)$, and this will be the
conclusion of the proof.

In the sequel we will denote by $h_i$ the real numbers $h_{M_i}
(P_i)$.

Take $d$ to be a very big and divisible positive integer; let
$d_i$ be integers such that $d_i h_i \sim d$ and such that
$$
{h_2 \over h_1} > {d_1 \over d_2}
$$
(in order to keep the proof as readable as possible we avoid to
introduce more small constants).

Let $f$ be the section whose existence is assured by theorem
\ref{smallsection}.

The hypotheses of theorem \ref{indextheorem} are satisfied
consequently the index of $f$ at $P_1 \times P_2$ is smaller then
$\epsilon_2$. Let $\tilde {\cal X} \rightarrow {\cal X}_1 \times
{\cal X}_2$ be the blow up of the ideal ${\cal
I}_{\underline{\vartheta} , \delta - \epsilon_2 , \underline d}$ and
$E_{\delta - \epsilon_2}$ (notations as in \S 5) be the exceptional
divisor; let $\tilde P : B \rightarrow \tilde{\cal X}$ be the strict
transform of $P := P_1 \times P_2$. We apply theorem
\ref{gencauchyineq} and deduce the existence of an absolute constant
$A_8$, a couple of indices $(i_1, i_2)$ and a non zero section
$\tilde f \in H^0 (\tilde P , M_1^{d_1} \otimes M_2^{d_2} \otimes
\omega_{{\cal X}_1 / B}^{i_1} \vert_{P_1} \otimes \omega_{{\cal X}_2
/ B}^{i_2} (-E_{\delta - \epsilon_2}))$
 such that ${i_1 \over d_1} \cdot \vartheta_1
+ {i_2 \over d_2} \cdot \vartheta_2 \leq \epsilon_2$ and
$\sup_{\sigma \in M_{\infty}} \{ \log \Vert \tilde f \Vert_{\sigma}
\} \leq A_8 \cdot T (D_1) \cdot T(D_2) (d_1 + d_2)$.

\label degreeofE. lemma\par\lemma Let $v\in M_K$. Then there is a
positive constant $C$ depending only on $v$ (and on the $\X_i$'s,
but independent on the $D_i$'s, $P_i$'s, $d_i$'s) and a couple
$(j_1^v;j_2^v)$ such that
${{j_1^v}\over{d_1}}\cdot\vartheta_1+{{j_2^v}\over{d_2}}\geq
\delta-\epsilon_2$ and
$$-\log\Vert E_{\delta-\epsilon_2}\Vert_v(\tilde P)\geq
j_1^v\lambda_{D_1,v}(P_1)+j_2^v\lambda_{D_2,v}(P_2)-C(d_1+d_2).$$
\endlemma

\Proof We prove the case when $v\in M_{\infty}$; when $v$ is
finite the proof is similar (and even easier).

We use the results and the notations of the proof of theorem
\ref{localcauchy}.

We may suppose that there is a couple $(j_1^v;j_2^v)$ with
${{j_1^v}\over{d_1}}\cdot\vartheta_1+{{j_2^v}\over{d_2}}\geq
\delta-\epsilon_2$ such that $\tilde P_v$ belongs to the open set
$\vert u_{j_1,j_2}\vert <(1+\epsilon)\vert u_{j^v_1,j^v_2}\vert$
for every $(j_1,j_2)$ in $H$. From the formula \ref{normofE}
$$\Vert E_{\delta-\epsilon_2} \Vert = \vert z_1^{j^v_1} \cdot z_2^{j^v_2}\vert  \cdot \sqrt{\sum_{(j_1 , j_2) \in H}
\left|{ {\vert u_{j_1 , j_2} \vert \cdot \rho_1^{j_1} \cdot
\rho_2^{j_2}}\over {\vert u_{i_1 , i_2} \vert}} \right|^2} \, .$$
Thus
$$\log\Vert E_{\delta-\epsilon_2}\Vert \leq j_1^v\log\vert z_1\vert+j_2^v\log\vert
z_2\vert+ C(d_1+d_2).$$ The conclusion follows from the properties
of the Weil functions.

\smallskip

From the lemma above, we deduce the existence of a constant $C$
such that
$$\widehat{\deg}(\tilde P^\ast(E_{\delta-\epsilon_2}))\geq
\sum_{v\in S}
j_1^v\lambda_{D_1,v}(P_1)+j_2^v\lambda_{D_2,v}(P_2)-C(d_1+d_2).$$
Consequently, from the hypotheses
$$\eqalign{\widehat{\deg}(\tilde P^\ast(E_{\delta-\epsilon_2}))&\geq \sum_{v\in S}
j_1^v\cdot\varphi(v)\cdot\vartheta_1\cdot
h_1+j_2^v\cdot\varphi(v)\cdot\vartheta_2\cdot d_2\cdot
h_2-C(d_1+d_2)\cr &=\sum_{v\in S}
\left({{j_1^v}\over{d_1}}\cdot\vartheta_1\cdot
d_1h_1+{{j_2^v}\over{d_2}}\cdot\vartheta_2\cdot
h_2\right)\cdot\varphi(v)-C(d_1+d_2)\cr &\geq
d(2+\epsilon_1-\epsilon_2)- C(d_1+d_2).\cr}$$

Thus we deduce
$$
\eqalign{ - \ &A_8 \cdot T(D_1) \cdot T(D_2) (d_1 + d_2) \cr \leq
\ &d_1 \cdot h_1 + d_2 \cdot h_2 + i_1 \cdot h_{\omega_{{\cal X}_1
/ B}} (P_1) + i_2 \cdot h_{\omega_{{\cal X}_2 / B}} (P_2) -
\widehat{\deg} (\tilde P^* (E_{\delta - \epsilon_2})) \cr \leq \
&((2 + 2 \epsilon_2 \cdot A_3) - (2 + \epsilon_1 - \epsilon_2))
\cdot d +C(d_1+d_2).\cr }
$$
From this and by our choice of the $\epsilon_i$'s, we deduce the
existence of a constant $A_9$ such that
$$
-A_9 \cdot T(D_1) \cdot T(D_2) \cdot \left( {1 \over h_1} + {1
\over h_2} \right) \leq -\epsilon_3
$$
where $\epsilon_3 = \epsilon_1 - (1 + 2A_3) \cdot \epsilon_2$; thus
$$
h_1 \leq {2 \cdot A_9 \over \epsilon_3} \cdot T (D_1) \cdot T
(D_2)
$$
and from this the conclusion follows.

\endsection

\

\

\centerline{\bf References}

\bigskip

\item{[Bi]} Bilu, Y. {\it Quantitative Siegel's theorem for Galois coverings}. Compositio Math. {\bf 106}
(1997), no. 2, 125--158.

\item{[B1]} Bombieri, E. {\it On the Thue-Siegel-Dyson theorem}.
Acta Math. {\bf 148} (1982), 255-296.

\item{[B2]} Bombieri, E. {\it The Mordell conjecture revisited}.
Ann. Scuola Norm. Sup. Pisa Cl. Sci. (4) {\bf 17} (1990), no. 4, 615-640.

\item{[B3]} Bombieri, E. {\it Effective Diophantine approximation on ${\bb G}_m$} .
 Ann. Scuola Norm. Sup. Pisa Cl. Sci. (4) {\bf 20} (1993), no. 1, 61-89.

\item{[BVV]} Bombieri, E., Van der Poorten, A.J., Vaaler, J.D.
{\it Effective measures of irrationality for cubic extensions of number fields}.
Ann. Scuola Norm. Sup. Pisa Cl. Sci. (4) {\bf 23} (1996), no. 2, 211-248.

\item{[BC]} Bombieri, E., Cohen, P.B. {\it Effective Diophantine approximation on ${\bb G}_M$. II}.
Ann. Scuola Norm. Sup. Pisa Cl. Sci. (4) {\bf 24} (1997), no. 2, 205-225.

\item{[Bo]} Bost, J.-B. {\it Algebraic leaves of algebraic foliations over number fields}.
Publ. Math. Inst. Hautes \'Etudes Sci. {\bf 93} (2001), 161-221.

\item{[BGS]} Bost, J.-B., Gillet, H., Soul\'e, C. {\it Heights of projective varieties and positive Green forms}.
J. Amer. Math. Soc. {\bf 7} (1994), no. 4, 903-1027.

\item{[CZ]} Corvaja, P. Zannier, U. {\it A subspace theorem approach to integral points on
curves}. C. R. Math. Acad. Sci. Paris {\bf 334} (2002), no. 4,
267--271.

\item{[D]} Deligne, P. {\it Le d\'eterminant de la cohomologie}. Current trends in arithmetical algebraic
geometry (Arcata, Calif., 1985),  Contemp. Math., 67, Amer. Math.
Soc., Providence, RI, 1987, 93--177.

\item{[Ev]} Evertse, J.-H. {\it An explicit version of Faltings' product theorem and an improvement of Roth's lemma}.
Acta Arith. {\bf 73} (1995), no. 3, 215-248.

\item{[Fa1]} Faltings, G. {\it Endlichkeitss\"atze f\"ur abelsche Variet\"aten ber Zahlk\"orpern}.
Invent. Math. {\bf 73} (1983), no. 3, 349-366.

\item{[Fa2]} Faltings, G. {\it Diophantine approximation on abelian varieties}.
Ann. of Math. (2) {\bf 133} (1991), no. 3, 549-576.

\item{[Ga]} Gasbarri, C. {\it  Some topics in Arakelov theory of arithmetic surfaces}.
Number theory, I (Rome, 1995). Rend. Sem. Mat. Univ. Politec. Torino
{\bf 53} (1995), no. 3, 309--323.

\item{[HS]} Hindry, M., Silverman, J.H. {\it Diophantine geometry. An introduction. Graduate
Texts in Mathematics}, {\bf 201}. Springer-Verlag, New York, 2000.

\item{[MB]} Moret-Bailly, Laurent. {\it M\'etriques permises}.
Seminar on arithmetic bundles: the Mordell conjecture (Paris, 1983/84). Ast\'erisque No. 127 (1985), 29-87.

\item{[Se]} Serre, J.P. {\it Lectures on Mordell Weil Theorem}.
Friedr. Vieweg \& Sohn, (1989).

\item{[Su]} Surroca, A. {\it Siegel's theorem and the $abc$ conjecture}.
Riv. Mat. Univ. Parma ({\bf 7}) 3* (2004), 323-332.

\item{[Sz]} Szpiro, L. {\it Pr\'esentation de la th\'eorie d'Arakelov}.
Current trends in arithmetical algebraic geometry (Arcata, Calif., 1985), 279-293,
Contemp. Math., {\bf 67}, Amer. Math. Soc., Providence, RI, 1987.

\item{[Vo1]} Vojta, P. {\it Dyson's lemma for products of two curves of arbitrary genus}.
Invent. Math. {\bf 98} (1989), no. 1, 107-113.

\item{[Vo2]} Vojta, P. {\it Siegel's theorem in the compact case}. Ann. of Math. (2) {\bf 133} (1991), no. 3, 509-548.

\item{[Vo3]} Vojta, P. {\it A generalization of theorems of Faltings and Thue-Siegel-Roth-Wirsing}.
J. Amer. Math. Soc. {\bf 5} (1992), no. 4, 763-804.

\

\item {Address:} C.~Gasbarri;   Dipartimento di Matematica
dell'Universit\`a di Roma ``Tor Vergata", Viale della Ricerca
Scientifica, 00133 Roma (I).

E--mail: gasbarri@mat.uniroma2.it

\bye